\documentclass[12pt]{article}
\usepackage[T1]{fontenc}
 \usepackage{stix}
 \usepackage{libertine,libertinust1math}
\usepackage[dvipsnames]{xcolor}
\usepackage{pdfpages}
\usepackage{xfrac}
\usepackage{sgame}
\usepackage{accents}
\usepackage{tikz-cd}
\usepackage{float}
\usepackage[round]{natbib}   
\usepackage{multirow}
\usepackage{dutchcal}
\usepackage{pdfpages}
\usepackage{sgame}
\usepackage{mathtools}
\usepackage{amsthm}
\usepackage{enumitem}
\usepackage[margin=1in]{geometry}
\usepackage{caption}
\usepackage{subcaption}
\usetikzlibrary{calc}
\usetikzlibrary{arrows}

\newcommand{\ubar}[1]{\underaccent{\bar}{#1}}

\newtheorem{theorem}{Theorem}[section]
\newtheorem{proposition}[theorem]{Proposition}
\newtheorem{lemma}[theorem]{Lemma}
\newtheorem{corollary}[theorem]{Corollary}
\newtheorem{claim}[theorem]{Claim}

\theoremstyle{definition}

\newtheorem{definition}[theorem]{Definition}
\newtheorem{example}[theorem]{Example}

\usepackage{hyperref}
\hypersetup{
    colorlinks=true,
    linkcolor=Violet,
    filecolor=Violet,      
    urlcolor=Violet,
    citecolor=Violet,
}


\DeclareMathOperator\supp{supp}
\begin{document}

\author{Pak Hung Au\thanks{%
Department of Economics, Hong Kong University of Science and Technology. Email: \href{mailto:aupakhung@ust.hk}{aupakhung@ust.hk}.} \and Mark Whitmeyer\thanks{Department of Economics, Arizona State University. Email: \href{mailto:mark.whitmeyer@gmail.com}{mark.whitmeyer@gmail.com}.} }
\title{Attraction versus Persuasion: Information Provision in Search Markets\thanks{
We thank the editor and five anonymous referees for many excellent comments. We are also grateful to Mark Armstrong, V. Bhaskar, Laura Doval, Ay\c{c}a Kaya, Willie Fuchs, Vasudha Jain, Elliot Lipnowski, Peter Norman, Larry Samuelson, Vasiliki Skreta, Isaac Sonin, Joseph Whitmeyer, Thomas Wiseman, Kun Zhang, and various audiences for their feedback. The excellent research assistance by Ke-cheng Hsu and Weixuan Zhou is gratefully acknowledged. This paper combines and subsumes the papers ``Competition in Designing Pandora's Boxes,'' and ``Persuading a Consumer to Visit;'' the latter of which comprised a chapter of the second author's dissertation at UT Austin. The first author gratefully acknowledges the financial support by the Lee Heng Fellowship. The second author thanks the University of Bonn, where he was supported under the DFG Project 390685813. Edited by Emir Kamenica.}}
\date{\today }
\maketitle

\begin{abstract}
We consider a model of oligopolistic competition in a market with search frictions, in which competing firms with products of unknown quality advertise how much information a consumer's visit will glean. In the unique symmetric equilibrium of this game, the countervailing incentives of attraction and persuasion yield a payoff function for each firm that is linear in the firm's realized effective value. If the expected quality of the products is sufficiently high (or competition is sufficiently fierce), this corresponds to full information--firms provide the first-best level of information. If not, this corresponds to information dispersion--firms randomize over signals.
\end{abstract}

\newpage

\section{Introduction}


Firms often have a great deal of control over how much information consumers' inspections will bring. Car sellers choose whether to offer test drives, software companies decide the length of free trial periods and the set of specific functions to include in promotional versions, newspapers limit how many free articles consumers can access, and book vendors specify how many and which pages consumers may sample for free. In choosing how much information to provide, a principal objective of each firm is clearly \textit{persuasion}: each wants consumers to select its product over those of its competitors.

Moreover, information about these inspections, such as the trial duration, functional restrictions, and number of free articles and/or pages available, are typically made available to prospective consumers \textit{before} they encounter the products. Firms advertise to consumers what information their probationary interactions with the product will bring, which affects their search behavior. Another central objective is therefore \textit{attraction}: each firm wants consumers to consider it (visit and inspect its product) early in their search.

In this paper, we investigate a series of fundamental questions: how are firms' information provision policies shaped by market competition? How do search frictions affect the intensity of competition? Would otherwise identical firms adopt common or idiosyncratic information policies? To answer these, we study a single-product oligopoly setting in which several firms compete by designing how much information a representative consumer obtains about their respective products through her (costly) inspections.

In our model, each \textit{ex ante} identical firm has a product of uncertain quality, which is either high or low.\footnote{Alternatively, this quality can be interpreted as the consumer's match value with the product. We use quality throughout in order to easily differentiate between it and the reservation value assigned to the firm in the consumer's search problem.} The quality distribution is independent and identically distributed across firms, with the probability of high quality denoted by $\mu $. Each firm has complete freedom over how much information about its product's quality a consumer's visit reveals: each firm simultaneously chooses and commits \textit{ex ante} to a signal or experiment, the outcome of which is revealed to the consumer only upon visiting that firm. These chosen signals are publicly posted and are observed by the consumer prior to commencement of her search.\footnote{\label{fn2}In justifying the relevance of their investigation of price-directed search, \cite{az} observe that ``price comparison websites are now a major part of the retailing landscape.'' Similar tools allow consumers to observe the information provision policies of firms before they inspect their products. Not only do firms advertise the length, breadth, and depth of their trial offers but there are also readily accessible lists online that make it easy to compare these disclosure policies.} The consumer needs at most one product.

Importantly, we abstract away from price competition and focus on the information provision problem of the firms. This is justified if the competing firms are dealers in some product whose price is set by a central office. Such firms can control the level of information that they provide but not the price. Similarly, platforms that host third party products often significantly limit the extent through which these third parties can compete via prices. Examples of this include the Apple App store, which has a fixed number of price tiers from which application developers can choose; and the Nintendo Switch eShop, which does not allow sellers to price titles below \$1.99.\footnote{An alternative rationale for our restriction to information competition is that the literature on consumer search has primarily explored the use of price as an instrument when information structures are fixed. We aim to understand the inverse: how do firms compete through information provision in a search environment (when prices are fixed)?} 

Knowing only the collection of signals posted but not their realizations, the consumer decides which firms to visit and in what order, at a search cost of $c>0$ per visit. As is standard in the consumer search literature, we assume that the consumer must visit a firm before she can buy from it and that recall is free: having visited a firm, the consumer may always return to select that firm's product. Consequently, at each stage in her sequential search, the consumer has two decisions to make: whether to stop (by selecting one product or the outside option) or continue, and whom to visit next if she continues.\footnote{Our model extends \cite{wei} by endogenizing the collection of prize distributions from which the consumer samples.}


As hinted above, there is an important tension present in the model that drives our results. Namely, when choosing its signal, a firm trades off between \textit{persuasion} and \textit{attraction}. Conditional on the consumer visiting, a firm wishes to maximize the chance that it is selected. The persuasion incentive encourages pooling of beliefs above the stopping threshold, as firms try to maximize the conversion rate from visits to purchases. On the other hand, committing to provide more information improves a firm's position in the consumer's search order. \textit{Ceteris paribus}, the earlier a firm is visited, the better for the firm, since ranking lower in the search order implies a greater chance that the consumer stops elsewhere before visiting it. In contrast to the persuasion incentive, the attraction incentive encourages the spreading of beliefs, in order to entice the consumer to visit. A firm's optimal signal is determined by the interplay of these two forces.

When the average quality, $\mu $, of the products or the number of competing firms, $n$, is sufficiently high, the attraction motive dominates, and the unique symmetric equilibrium is one in which each firm chooses a fully informative signal.\footnote{This result rings true, in the sense that there are many markets in which firms provide a large amount of information to visiting consumers. The fine piano purveyor, Steinway \& Sons, has practice rooms in its showrooms, many car dealerships allow test-drives and some even allow prospective buyers to keep cars overnight, upscale clothing stores include changing rooms (and mirrors) for trying on their products, and anyone who has wandered through the perfume section of a department store can attest that the perfume sellers provide ample olfactory evidence about their wares.} There is no profitable deviation from full information, as any other signal would ensure that the deviating firm be visited last, a rare event due to the surfeit of high quality products in the market. Conversely, if $\mu $ or $n$ is not high, there are no symmetric equilibria in pure strategies. Because the average quality is low and there are not many competitors, the persuasion motive is more important. A firm can deviate profitably from providing full information by being uninformative: even though it will be visited last--indeed, it will only be visited if every other firm is low quality--it will be selected for sure if visited. Nevertheless, the attraction incentive still remains and precludes the existence of any other symmetric pure strategy equilibrium, as firms can always deviate profitably by providing slightly more information and moving to the top of the consumer's search order.

In order to characterize the unique symmetric equilibrium outside of the high average quality case, we use recent results from \cite{am} and \cite{choi2}, who show how the sequential problem of \cite{wei} can be reformulated as a static discrete choice problem in which the consumer selects the firm with the highest realized \textit{effective} value. Adapting the concavification method of \cite{kam} to this setting, we establish that there is a unique symmetric equilibrium distribution over effective values, which necessarily begets a payoff function for each firm
that is linear in the firm's realized effective value. Importantly, this distribution requires firms to randomize over their choices of signals, i.e., our model generates \textit{information dispersion}.\footnote{Note that by information dispersion, we mean the informational analog of price dispersion, \textit{viz.}, the variation in levels of information provision among identical firms in the same market.} 

Our prediction of information dispersion is novel and suggests potentially profitable avenues for empirical study. For instance, in the market for antivirus/security software, the major companies offer potential consumers free trials, but differ in the length of these previews: 7 days for AVG and Avast, 14 days for MalwareBytes, and 30 days for Norton, McAfee and Kaspersky. Similarly, there is variation in the length of free trial periods offered by music streaming services: 1 month for Google Play Music and Tidal, 2 months for Pandora, and 3 months for Spotify and Apple Music. The market for grammar-checking software features not only different free-trial periods, but also different functionalities (including word limits and style settings) in firms' probationary offers.\footnote{In each of these industries (antivirus software, music streaming services, and grammar-checking software), firms do the sort of advertising of disclosure policies noted in Footnote \ref{fn2}. Furthermore, adjusting these trial versions may be difficult or impossible in the short run, thus providing support for our simultaneous-move specification.}%

To highlight the significance of the attraction motive in driving our results, we also investigate two benchmarks in which only the persuasion incentive is present, albeit for different reasons. In Section \ref{costlessbench}, we explore the case in which the consumer is able to observe the signal realizations of all of the firms for free (corresponding to setting $c=0$ in our main model). In the absence of the attraction motive, firms have less incentive to provide information; and in the unique symmetric equilibrium, firms do not provide full information. Consequently, our results imply that when it is important to be visited early ($\mu $ or $n$ is high), consumers can benefit from a small positive search cost, which engenders the attraction incentive.

In Section \ref{diamond}, we consider another benchmark in which the
consumer can observe a firm's choice of signal only after paying it a visit.
When information is hidden in this manner, the attraction motive is
completely absent, and we point out a dramatic \textquotedblleft
informational Diamond paradox.\textquotedblright\ Namely, the only
equilibrium outcome is that firms provide no useful information and the
consumer does not actively search. Deviations to other signals cannot be
observed in advance and so the consumer's search order is determined
entirely by her conjectures and not the actual signals. The pooling
incentive is all that remains, which eliminates any purported equilibria
with active search. Again, there is no information dispersion as all firms
provide the monopoly level of information.

For tractability, our main model assumes \textit{ex ante} homogeneous firms, so it is natural to focus on symmetric equilibria. For simplicity, we also assume the consumer's outside option is irrelevant. In Section \ref{extensions}, we discuss the consequences of relaxing these restrictions. Specifically, we illustrate that asymmetric equilibria are possible in some parameter regions, analyze competition between two asymmetric firms, show that the introduction of a non-negligible outside option leaves the qualitative features of the equilibrium intact, and argue that the principal findings of the main setting extend to an environment in which each firm's quality is distributed continuously on an interval.

This section concludes with a brief discussion of related work. The model is set up in Section \ref{themodel}. Section \ref{prelim} reports some preliminary observations and explains how the game under study can be reformulated into a more tractable one. Results on equilibrium existence, uniqueness and characterization are detailed in Section \ref{eqsy}. Section \ref{bench} illustrates the major economic forces at work by considering two benchmark models. A number of extensions are considered in Section \ref{extensions} before Section \ref{discuss} concludes. All proofs are left to the appendices, unless stated otherwise.

\subsection{Related Work}

\label{relwork}

There are still relatively few papers that explore information design in search settings. \cite{board} consider a setting in which sellers compete by designing experiments and buyers search sequentially. In contrast to our paper, the sellers' experiments are not publicly posted and do not direct the buyers' search, eliminating the tension between attraction and persuasion. They show that under certain conditions, the monopoly outcome (no active search) is a unique equilibrium. As in our hidden information benchmark, firms want to pool information just above a consumer's stopping threshold. 

\cite{hu} explore consumer-optimal information structures in the sequential (undirected) search framework of \cite{asher}. They find that consumer
welfare is maximized by a signal that generates a (conditional) unit-elastic demand. A signal that generates unit-elastic demand is also optimal in \cite{pease}, who look at consumer-optimal signals in a monopoly problem for search goods, and in other papers--\cite{Roes}, \cite{Condor}, and \cite{yang}--that study information structures in bilateral trade.

There are also several papers that explore competition through information provision when there are no search frictions ($c=0$). \cite{spiegler} explores a closely related (mathematically) scenario, \cite{cotton} and \cite{al} derive results that characterize the two-player solution to this problem of frictionless competitive information provision, and \cite{Au2} characterize the unique (symmetric) equilibrium in the $n$-player game. \cite{Koessler} look at a general setting in which multiple persuaders provide information about their own dimension of some multidimensional state.

A natural point of comparison for this paper is the collection of papers that explore price-directed search.\footnote{This list includes \cite{az}, \cite{choi2}, \cite{ding}, and \cite{haan}.} In each such paper, there is a tension inherent to firms' pricing decisions. Setting a lower price makes a firm more likely to be visited early as well as make a sale if visited, yet lowers the firm's profit directly. As search frictions increase, it becomes more important to attract (and retain) consumers, which drives prices down. At first glance, this trade-off seems like a direct analog of the persuasion/attraction conflict in this paper. However, there are important differences. In particular, lower prices help with both persuasion and attraction--consumers are both more likely to visit \textit{and} to stop as a firm's price shrinks. The idea that firms can increase the chance of being selected if visited, at the expense of being visited in the first place, is completely absent from this literature that focuses on pricing; yet is the driving feature in our model. To put
differently, in both price and information-directed competition in consumer search, the incentive to be visited early is fundamental, yet the costs of
such prominence are different.


Also related are the works that explore other varieties of frictions in markets. \cite{wol} allow firms to increase the length of time it takes for consumers to learn their prices, to the consumers' detriment. Just as search frictions can improve consumer welfare in our model, the rational inattention literature points out that so too can informational frictions. \cite{eliaz} unearth an incentive for prominent firms to lower prices and therefore stay \textquotedblleft under consumers' radar,\textquotedblright\ in contrast to our firms, who compete to get ``on consumers' radar.'' \cite{ravid} finds that attention costs strictly benefit a buyer in a bargaining setup, who obtains bargaining power as a result of his inattention.

\section{The Model}

\label{themodel}

\label{model}

There is one consumer and $n$ \textit{ex ante} identical single-product firms indexed by $i$. Each firm's product has an uncertain quality (or match value) to the consumer of either $0$ or $1$. These qualities are identically and independently distributed, with $\mu \in \left( 0,1\right) $ being the prior probability that the quality is $1$. The consumer has unit demand, and her \textit{ex post} payoff of consumption is normalized to the quality of the product consumed. For simplicity, we assume that the consumer has no outside option.\footnote{In Section \ref{extensions}, we briefly discuss the impact of a positive outside option.} A firm receives a payoff normalized to $1$ if the consumer picks its product and $0$ otherwise. All players are risk neutral.

At the beginning of the game, neither the consumer nor the firms know the quality realizations of the firms' products. Each firm simultaneously commits to a signal, $\pi_{i}\colon \left\{0,1\right\} \rightarrow \Delta (S)$ with some space of signal realizations $S$. The primary focus of this paper is the scenario in which the chosen signals are publicly posted and therefore shape the consumer's behavior directly.\footnote{In Section \ref{diamond} we study a setting in which the chosen signals are hidden and can only be observed after incurring the search cost. Naturally, in this case, the consumer conjectures firms' signal choices, which must be correct at equilibrium.} Our analysis allows firms to play mixed strategies (i.e., randomize over signals). In this case, the firms' randomizations are resolved simultaneously before the consumer embarks on her sequential search.

Guided by the informativeness of the posted signals, the consumer learns about the firms' product qualities by visiting the firms and observing their signal realizations in sequence. Each such inspection requires a search cost of $c>0$. After observing the signal realization of a firm, the consumer updates her prior to the posterior expected quality of the product offered by that firm. At any stage of her sequential search, the consumer can stop her search by buying from any previously visited firm; recall is assumed to be free. Alternatively, she can continue her search by visiting more firms or stop her search by collecting the outside option. We assume that $c\leq \mu$, so that search is not strictly dominated for the consumer.

The solution concept is subgame-perfect Nash equilibrium. Given the \textit{ex ante} symmetry of the firms, it is natural to focus on symmetric equilibria in which (i) all firms adopt a common (possibly mixed) strategy, and (ii) the consumer uses a tie-breaking rule that treats all firms identically.

\section{Preliminary Analysis}

\label{prelim}

The purpose of this section is to simplify the game set up in the last
section. First, we note that a firm's strategy space can be expressed as the set of distributions over posterior qualities and that the consumer's sequentially optimal search strategy takes a simple form. We use a stylized example to compare and contrast the effects of attraction and persuasion on a firm's information provision. Second, we explain the necessity of considering mixed strategies by noting that there are no pure-strategy equilibria in a large region of the parameter space. Finally, using the fact that the consumer's decision can be expressed as a static discrete choice problem, we show how the strategic interaction between the firms can be condensed and simplified.

\subsection{Basics}

This subsection reformulates the game set up in Section \ref{themodel} into
one of competition in the design of distributions over posteriors. Instead
of choosing signals directly, the strategy space of each firm can be
redefined without loss to be the set of feasible distributions over
posterior (expected) qualities. To that end, we introduce the following
definition.

\begin{definition}
A distribution over posteriors, $F_{i}$, is feasible if it has mean $\mu $ and its support is a subset of $\left[ 0,1\right] $. $\mathcal{F}$ denotes the set of feasible distributions.
\end{definition}

As a signal realization affects payoffs only through its induced (posterior) expected product quality, it is a standard result in the literature that we
may define a firm's pure strategy to be a feasible distribution over posterior qualities. Consequently, a generic mixed strategy of a firm is a randomization over the set of feasible distributions, and the set of mixed strategies is $\Delta \left( \mathcal{F}\right) $. Recall that if a firm plays a mixed strategy, its pure strategy realization occurs and becomes public \textit{before} the consumer begins her search.

Next, given the $n$ posted distributions over posterior qualities, the sequentially optimal selection and stopping rule is identified by \cite{wei}. A brief recap of his finding is useful. For any distribution $F_{i}$ chosen by firm $i$, define the corresponding reservation value, $U\left(F_{i}\right) $, implicitly as the solution to the following equation (in $U$): 
\[c=\int_{U}^{1}\left(x-U\right) dF_{i}\left(x\right) \text{ .}
\label{reservation value} \tag{$1$}\]
The set of feasible reservation values $\left\{ U\left( F_{i}\right)\right\} _{F_{i}\in \mathcal{F}}$ is bounded between $\ubar{U} \equiv \mu -c$ and $\bar{U}\equiv 1- c/\mu$. The lower bound is induced by any feasible distribution whose support is entirely (weakly) above $\ubar{U}$, one of which is the degenerate distribution at $\mu$ (which corresponds to a completely uninformative signal). The upper bound is uniquely induced by the feasible distribution supported on $\left\{0,1\right\} $ (which corresponds to a fully revealing signal). It is not difficult to see that any intermediate reservation value can be achieved by some feasible distribution,\footnote{For instance, by some convex combination of the two distributions above.} so $\left\{ U\left( F_{i}\right)\colon F_{i}\in \mathcal{F}\right\} =\left[ \ubar{U},\bar{U}\right] $. Moreover, the reservation value rewards informativeness: for any pair of feasible distributions $F_{i}$ and $G_{i}$, if $F_{i}$ is a mean-preserving spread of $G_{i}$, then $U\left( G_{i}\right) \leq U\left( F_{i}\right) $.\footnote{This is shown, e.g., in Corollary 20 of \cite{kohn1974theory}.}

The optimal strategy of the consumer (\hypertarget{pandora}{Pandora's rule}) is as follows.
\begin{itemize}[noitemsep,topsep=0pt]
\item Selection rule: If a firm is to be visited and examined, it should be the unvisited firm with the highest reservation value.
\item Stopping rule: Search should be stopped whenever the maximum reservation value of the unvisited firms is lower than the maximum sampled
reward or the outside option.
\end{itemize}

As we focus on symmetric equilibria, we select the (sequentially) optimal strategy in which the consumer adopts a fair tie-breaking rule throughout the search process. With the consumer's search behavior pinned down, we can we can now focus our attention on the strategic interaction between the firms.

\subsubsection{Attraction versus Persuasion}\label{311}

This subsection uses a simplified setting that differs slightly from the main model outlined in Section 2 to illustrate and disentangle the attraction and persuasion incentives faced by a firm. Let the market consist of just two firms. In order to focus on the information choice of firm 1, we assume that firm 2 is non-strategic and chooses a distribution over posteriors with support $\left\{0, p_2\right\} $ for some $p_2 > \mu $. This yields the reservation value $p_2\bar{U}$ ($< p_2$). We also assume that the consumer favors firm 1 whenever she is indifferent.

To remove the attraction motive, suppose first that firm 1 is prominent and free to visit, so that the consumer always visits it first. Firm 1's objective is simply to persuade the consumer to buy from it rather than its competitor. If firm 1's posterior quality exceeds $p_2\bar{U}$, the consumer buys from it right away without visiting firm 2. Otherwise, the consumer also visits firm 2 but firm 1 still gets her business if firm 2's posterior quality is $0$. firm 1's payoff as a function of its realized posterior $p$ is
\[V\left(p\right) =
\begin{cases}
1-\frac{\mu }{p_{2}}, \quad &\text{if} \quad  p \in \left[0,p_{2}\bar{U}\right) \\ 
1,\quad &\text{if} \quad p\geq p_{2}\bar{U}
\end{cases} \text{ .}
\]
Following \cite{kam}, firm 1's
optimal distribution of posteriors has support $\left\{0,p_{2}\bar{U}\right\} $ if $p_{2}\bar{U}>\mu $ and is degenerate at $\mu$ (no information) otherwise.

How does the need to attract the consumer affect firm 1's signal design? Suppose firm 1 is no longer prominent so that the consumer's (optimal) shopping strategy follows \hyperlink{pandora}{Pandora's rule}. If firm 1 still adopted the distribution of posteriors above, it would no longer be visited first (as its reservation value would be below that of firm 2). Being visited last, firm 1 would make a sale only if firm 2's posterior quality was $0$, yielding it an expected payoff of $1-\mu/p_{2}$. To improve its payoff, firm 1 must promise more information so that the consumer is willing to visit it first. Mimicking firm 2's signal is one such increase, yielding a reservation value equal to that of its competitor. Being the first to be visited (recall that ties are broken in firm 1's favor), firm 1 makes a sale if it has posterior realization $p_2 \bar{U}$ or if both firms have a zero posterior realization. Consequently, its expected payoff improves to $\mu /p_{2}+\left(
1-\mu /p_{2}\right)^{2}>1-\mu /p_{2}$.\footnote{%
Example \ref{opt eff dist}, \textit{infra}, shows that this distribution is optimal for firm 1.}

The takeaway message from this investigation is as follows. For an individual firm, effective persuasion typically calls for partial information disclosure: the goal is to maximize the likelihood of realizing posteriors that are just good enough for persuasion to succeed. The need to attract the consumer's visit calls for information revelation beyond that that is required for effective persuasion alone. As the subsequent analysis shows, the interplay of the persuasion and attraction incentives plays the primary role in shaping the information policies in a competitive equilibrium.

\subsection{The Full Information Equilibrium}
\label{fullinfoeq}

This subsection focuses on pure-strategy equilibria and identifies the
conditions under which one exists as well as the form it takes. We begin
with a simple observation:

\begin{lemma}
\label{no partial pure eqm}There exist no symmetric pure strategy equilibria
in which any reservation value $U<\bar{U}$ is induced.
\end{lemma}

The rationale mirrors that of the Bertrand model of homogeneous goods in
which pricing above marginal cost cannot be sustained in equilibrium. Just
as a firm can ``undercut'' its rivals' marked-up prices to obtain a discrete jump up in its demand; here, a firm can ``overcut'' its rivals' less-than-full provision of information. More specifically, a firm can perturb its distribution slightly, choosing a nearly identical--but slightly more informative--distribution.\footnote{One way of doing this is by increasing the probability of values $0$ and $1$ by an arbitrarily small amount and decreasing the probability of a commensurate measure of interior values.} Doing so grants the firm a considerable edge on its competition (since it will move to the top of the consumer's search order), at a negligible loss of persuasion effectiveness (since its distribution is nearly unchanged). This results in a discrete gain in expected profit, so a symmetric partial-information equilibrium cannot be sustained. This leaves full information as the only candidate pure-strategy equilibrium.

\begin{proposition}
\label{2states} Define $\Bar{\mu}\equiv 1-\left( 1/n\right)^{\frac{1}{n-1}}$%
. A symmetric equilibrium in pure strategies exists if and only if $\mu \geq 
\Bar{\mu}$; i.e., the average quality, $\mu $, is sufficiently high or the
number of firms, $n$, is sufficiently large. In this equilibrium, all firms
provide full information.
\end{proposition}

This proposition details precisely the conditions under which the attraction incentive dominates the persuasion incentive for the firms. A high $\mu $ implies persuasion is likely to succeed, making it paramount for a firm to entice the consumer into visiting it--a failure to do so means the consumer is likely to stop her search at one of the firm's rivals before ever reaching it. A similar effect is at work if the number of rival firms is large--a low rank in the consumer's search order means a firm is unlikely to ever be visited, let alone make a sale. Consequently, in these cases, the attraction incentive dominates. If rival firms reveal full information, a departure guarantees the deviator the bottom rank in the consumer's search order, a devastating outcome when the average quality is high or there are a large number of competitors.

The relative importance of attraction is smaller when $\mu $ or $n$ is relatively low. Even if all of its rivals provide full information, an individual firm may find it profitable to provide less information: despite being ranked last, the firm still has a decent chance of eventually being visited, in which case it makes a sale with probability one. When the attraction motive no longer dominates, both the attraction and persuasion incentives play a part in shaping the equilibrium outcome. By Lemma \ref{no partial pure eqm}, the equilibrium involves mixed strategies, making its characterization more involved.

Interestingly, the existence of the full-information equilibrium does not depend on the consumer's search cost $c$ (as long as $0 < c<\mu$). This is because a firm's profitability of deviating from full information is independent of the search cost. In fact, if all firms but one provide full information, the optimal search strategy implied by \hyperlink{pandora}{Pandora's rule} is itself independent of the search cost: the consumer simply visits the full-information firms first, and stops if and only if she obtains a high quality realization.

In the next subsection, we make use of a remarkable discovery of \cite{choi2} and \cite{am}, who show that the consumer's sequential search can be formulated as a (static) discrete choice problem, which enables us to characterize and establish the uniqueness of the symmetric equilibrium in mixed strategies in a tractable way.\footnote{As it turns out, the full disclosure equilibrium exposed in this subsection corresponds to a special case of a more general result that we derive in the next section. The main result there, Theorem \ref{summary no u0}, establishes that the (essentially) unique equilibrium of the game (for all parameter values) must take a particular linear form, of which the full disclosure equilibrium serves as an avatar when $\mu $ is sufficiently high.}

\subsection{Reformulating the Game and Main Analysis}\label{secreform}

\label{reformulationsection}
Using the observation that the consumer's optimal shopping strategy affects firms' payoffs only through her eventual purchase decision, we model the strategic interaction between firms as competition in choosing distributions of effective values. The effective value of firm $i$ is defined as $W_{i}\equiv \min \left\{ p_{i},U\left( F_{i}\right) \right\} $, where $p_{i}\in \left[ 0,1\right] $ is the realized posterior quality and $F_{i}\in \mathcal{F}$ is the distribution over posteriors offered by firm $i$. \cite{choi2} observe that under Pandora's rule, there are two (and only two) scenarios in which firm $i$ makes a sale. First, if firm $i$'s effective value is its reservation value, the consumer purchases from it whenever she visits it, an event that occurs if and only if all of the other firms she visited earlier (due to their higher reservation values) realize posterior qualities below $U\left(F_{i}\right)$. Second, if firm $i$'s effective value is its posterior quality $p_{i}$, it makes a sale if and only if firms whose posterior qualities are greater than $p_{i}$ are not visited (due to their reservation values falling short of $p_{i}$). 

Putting these together yields the winning condition for firm $i$: its effective value realization, $\min \left\{ p_{i},U\left( F_{i}\right) \right\}$, beats the posterior qualities of those competitors visited before it (the former scenario) and the reservation values of those with higher posterior qualities (the latter scenario). That is, firm $i$ wins the consumer's business only if it has the highest effective value realization, i.e., $W_{i}=\max_{j=1,2,...,n}W_{j}$. Consequently, two (possibly mixed) strategies of a firm are payoff-equivalent as long as their implied effective-value distributions are identical, because they imply the same \textit{ex ante} probability of making the sale. This observation allows us to
reformulate the problem facing an individual firm as one where it chooses an effective-value distribution in order to maximize the probability that its realized effective value is greater than those of its competitors (and the consumer's outside option).

This transformation requires us to characterize the set of effective-value distributions that can be induced by a firm's strategy. Formally, we say an effective-value distribution is inducible if and only if there is a strategy, either pure or mixed, that can induce it. Consider first the simpler case of pure strategies. After a feasible distribution $F_{i}\in \mathcal{F}$ is chosen, but before the posterior realizes, firm $i$'s effective value, $W_{i}$, is a random variable with distribution
\[H\left(w; F_{i}\right) \equiv \begin{cases}
F_{i}\left( w\right) \quad &\text{if} \quad w < U\left(F_{i}\right) \\ 
1 \quad &\text{if} \quad w \geq U\left( F_{i}\right)
\end{cases}\text{ .}  \label{eff-value dist.}\tag{$2$}
\]
Effective-value distributions induced by pure strategies have several
notable features.

\begin{lemma}
\label{effective-value mean}Let $F_{i}$ be a feasible distribution over posteriors with reservation value $U\in \left[ \ubar{U},\bar{U}\right] $. Its induced effective-value distribution $H\left(\cdot ;F_{i}\right) $ has the following properties.
\begin{enumerate}[label={(\roman*)},noitemsep,topsep=0pt]
    \item \label{1stpart} $H\left( \cdot ;F_{i}\right) $ has an atom of size at least $c/\left(1-U\right)$ on the reservation value $U$. Moreover, $H\left(U;F_{i}\right) =1$.
    \item \label{2ndpart} $H\left( \cdot ;F_{i}\right) $ has mean $\ubar{U}=\mu -c$.
\end{enumerate}
\end{lemma}
The properties of $H\left( \cdot ;F_{i}\right) $ reported in the lemma above are straightforward consequences of the definition of an effective-value distribution (\ref{eff-value dist.}), the reservation-value equation (\ref{reservation value}) and the feasibility of $F_{i}$.

Let us briefly take a step back from the effective value reformulation and consider the ramifications of allowing firms to randomize over information structures. While allowing for mixed strategies does not expand the set of feasible distributions over posteriors (from an \textit{ex ante} point of view), it does--given the simultaneous nature of the signal posting--substantively expand the firms' (payoff-relevant) strategic capacities. First, it convexifies the set of effective-value distributions that can be induced deterministically, which is not, itself, convex.\footnote{Example 1.1 in the \href{https://whitmeyerhome.files.wordpress.com/2022/05/attraction_versus_persuasion_supplementary_appendix_jpe_final.pdf}{Supplementary Appendix} demonstrates this.} Second, a deterministic distribution over posteriors produces an atom at the reservation value in the induced effective-value distribution, making a firm's payoff highly sensitive to whether other firms' reservation values are above or below its reservation value. Mixed strategies absolve the firm from this requirement and hence this payoff sensitivity.

The following lemma provides an exact characterization of an inducible effective-value distribution by representing it using a reservation-value distribution and a collection of feasible distributions over posteriors, one for each relevant reservation value.

\begin{lemma}
\label{inducibility}An effective-value distribution $H_{i}$ is inducible if
and only if there is a reservation-value distribution $G_{i}\in \Delta
\left( \left[ \ubar{U},\bar{U}\right] \right) $, where, for each $U \in \supp\left\{ G_{i}\right\} $ there is a feasible distribution over
posteriors $F_{i}\left( \cdot ; U\right) \in \mathcal{F}$ with reservation
value $U$ such that for each $w\in \left[ 0,\bar{U}\right] $, 
\[
H_{i}\left( w\right) =G_{i}\left( w\right) +\int_{supp\left( G_{i}\right)
\cap (w,\bar{U}]}F_{i}\left( w; U\right) dG_{i}\left(U\right) \text{ .}
\label{effective value dist}\tag{$3$}
\]
\end{lemma}

The effective-value distribution formula (\ref{effective value dist}) can be
understood as follows. Take some $w\in \left[ 0,\bar{U}\right] $ and note
there are two scenarios in which firm $i$ realizes an effective value below $%
w$: either its reservation value realization (after the strategic mixing is resolved) falls below $w$, an event that happens with probability $G_{i}\left( w\right) $; or its reservation value exceeds $w$ but its posterior quality realization falls short of $w$. The probability of the latter scenario is captured by the integral in Equation (\ref{effective value dist}).

Part \ref{2ndpart} of Lemma \ref{effective-value mean} reveals that an inducible effective-value distribution has mean $\ubar{U}=\mu -c$.\footnote{This can be obtained by applying the Law of Iterated Expectations.} This mean condition; however, is not sufficient for inducibility. For instance, the only distribution over effective values with $\bar{U}$ in support is the binary distribution with support $\left\{0,\bar{U}\right\}$ (corresponding to full information).\footnote{Example 1.2 in the \href{https://whitmeyerhome.files.wordpress.com/2022/05/attraction_versus_persuasion_supplementary_appendix_jpe_final.pdf}{Supplementary Appendix} elaborates on this point.} Consequently, while it is tempting, we may not treat a firm's problem as a standard persuasion problem in which it just chooses a distribution over effective values with a specific mean condition. The next subsection studies a graphical approach to a firm's optimization problem in which it chooses from the set of inducible effective-value
distributions.

\subsubsection{Finding the Optimal Inducible Effective-Value Distribution by
Concavification}\label{modapproach}

In this subsection, we show how to adapt the concavification approach of \cite{kam} to this environment in order to find a firm's optimal inducible\footnote{Henceforth, to save space we omit this modifier.} effective-value distribution. To this end, we temporarily focus on the problem of a single firm and specify its payoff to be some function of its realized effective value $\Pi\colon \left[ 0,\bar{U}\right] \rightarrow \mathbb{R}$. For expositional simplicity, we assume that $\Pi$ is continuous.\footnote{Without the continuity assumption, an optimal distribution may not exist, but it is clear that the approach below still identifies distributions that deliver payoffs arbitrarily close to the supremum.}

The quest for the optimal effective-value distribution can be decomposed into two steps. We first identify the optimal effective-value distribution for each reservation value $U\in \left[\ubar{U},\bar{U}\right] $. Comparing the expected payoffs associated with each reservation value subsequently determines the overall optimum. The economics behind this approach is intuitive. A firm's choice of reservation value, $U$, caters directly to the attraction motive. In the first step, when $U$ is fixed, the attraction motive is completely absent and the firm is only concerned with persuading effectively. This persuasion problem differs from the standard one due to the constraints engendered by the firm's choice of $U$, as specified in Lemma \ref{effective-value mean}. While a high reservation value $U$ facilitates attraction, part \ref{1stpart} of the lemma implies that a high $U$ impairs the firm's persuasion ability by limiting the mass that it can allocate to effective values other than $U$. The second step, in which the firm chooses $U$, resolves this trade-off between attraction and persuasion.

\begin{figure}[tbp]
\centering
\includegraphics[scale=.15]{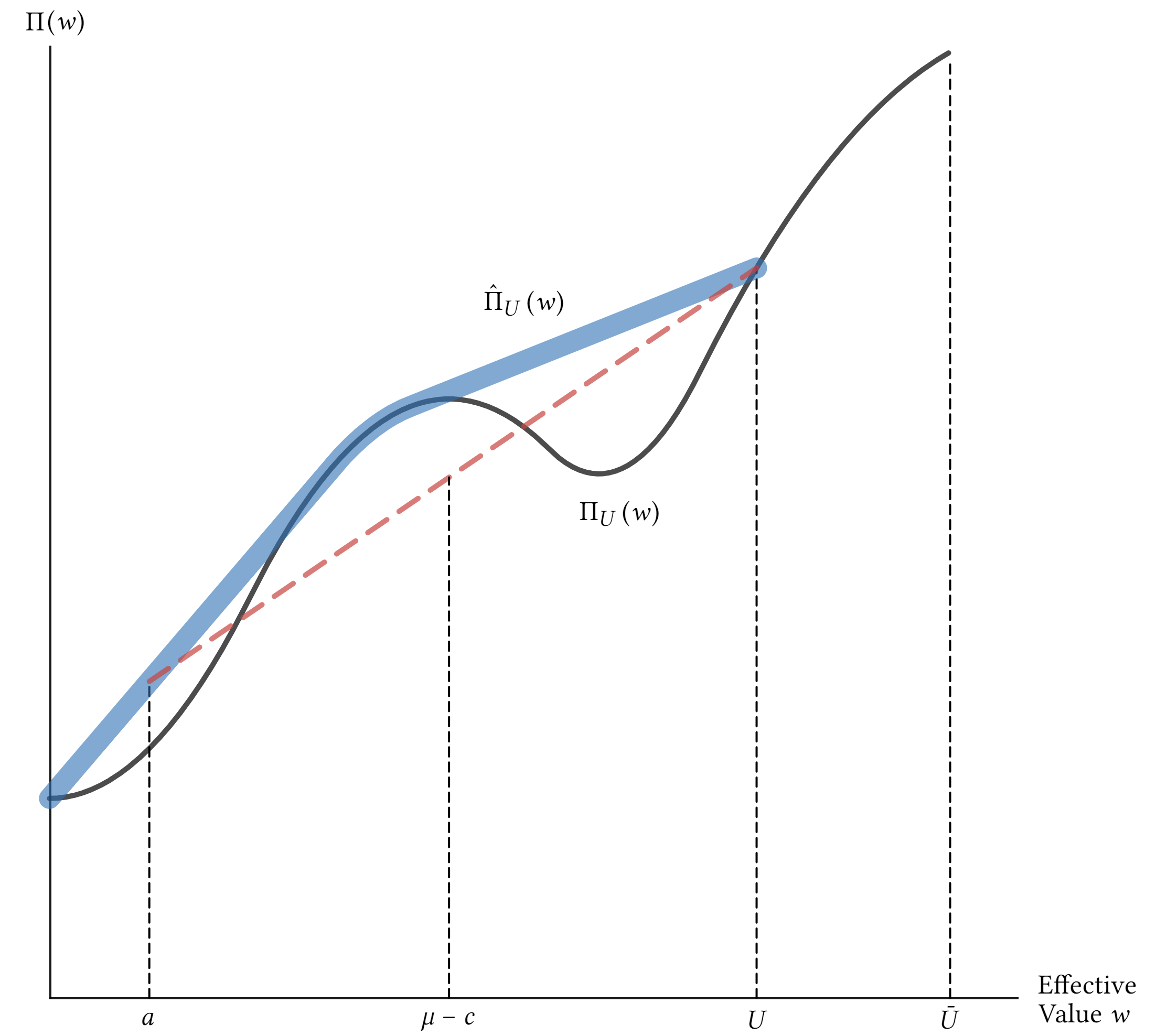}
\caption{An illustration of $\Pi $; its restriction to $\left[ 0,U\right] $, 
$\Pi _{U}$; and the concavification of $\Pi _{U}$, $\hat{\Pi}_{U}$.}
\label{fig1}
\end{figure}

We begin by fixing a reservation value $U\in \left[\ubar{U},\bar{U}\right]$ and solve for the effective-value distribution that optimizes persuasion. Part \ref{1stpart} of Lemma \ref{effective-value mean} states that a mass of at least $c/\left( 1-U\right) $ must be assigned to effective value $U$, which can be interpreted as a firm's fixed persuasion budget. Part \ref{2ndpart} of the lemma states that a firm's distribution must have a mean of $\ubar{U}$. Thus, a firm must also choose how to distribute the residual mass of $1-c/\left( 1-U\right) $ over the interval $\left[0,U\right] $. The conditional mean of this residual mass, or the flexible persuasion budget, $a\left( U\right) $, is pinned down by
\[\frac{c}{1-U}U+\left( 1-\frac{c}{1-U}\right) a\left( U\right) = \ubar{U} \quad \Leftrightarrow \quad a\left( U\right) =\frac{\left(1-U\right) \ubar{U}-cU}{1-U-c} \text{ .}  \label{baradef}\tag{$4$}\]
The optimal allocation of this flexible persuasion budget can be obtained by applying the concavification approach of \cite{kam} to the restriction of payoff function $\Pi $ to domain $\left[ 0,U\right] $ with a conditional mean set at $a\left( U\right) $. The (conditional) contribution of the flexible persuasion budget to the firm's payoff is $\hat{\Pi}_{U}\left(a\left( U\right) \right) $. Accordingly, the firm's payoff is the weighted sum of the contributions from the fixed and flexible persuasion budgets:
\[V\left( U\right) = \underbrace{\left( 1-\frac{c}{1-U}\right) \times \hat{\Pi}_{U}\left(
a\left( U\right) \right) }_{\text{contribution by flexible persuasion budget}}+\underbrace{\frac{c}{1-U}\times \Pi \left( U\right)}_{\text{contribution by fixed persuasion
budget}}\text{ .}\]
On the graph of $\Pi $, this value is the intersection of the vertical line $%
w=\mu -c$ and the line connecting $\left( a,\hat{\Pi}_{U}\left( a\right)
\right) $ and $\left( U,\Pi \left( U\right) \right) $. See Figure \ref{fig1}
for an illustration.

Having identified the optimal effective-value distribution for each reservation value, the overall optimum can be found by optimizing $V\left(U\right) $ over $U\in \left[\ubar{U},\bar{U}\right] $. Equation (\ref{baradef}) indicates that $a\left( U\right) $ varies inversely with $U$. The trade-off in the choice of $U$ is thus attraction versus persuasion: although a higher $U$ facilitates attraction, it leaves a smaller flexible persuasion budget $a\left( U\right) $, thereby impeding persuasion. The following proposition summarizes our graphical approach for effective-value distribution optimization.

\begin{proposition}
\label{optimal prize distribution} Denote $U^{\ast }\equiv \arg \max_{U\in \left[\ubar{U},\bar{U}\right] }V\left( U\right) $. The firm's optimal payoff is $V\left( U^{\ast }\right) $, which can be achieved by an effective-value distribution assigning mass $c/\left( 1-U^{\ast }\right)$ to $U^{\ast }$ and the remaining mass to values in the interval $\left[0,U^{\ast }\right] $ according to the construction of the concavification $\hat{\Pi}_{U^{\ast }}$ at $a\left( U^{\ast }\right) $.
\end{proposition}

We conclude this section with an example that illustrates the use of the
graphical approach.

\begin{figure}
    \centering
    \includegraphics[scale=.15]{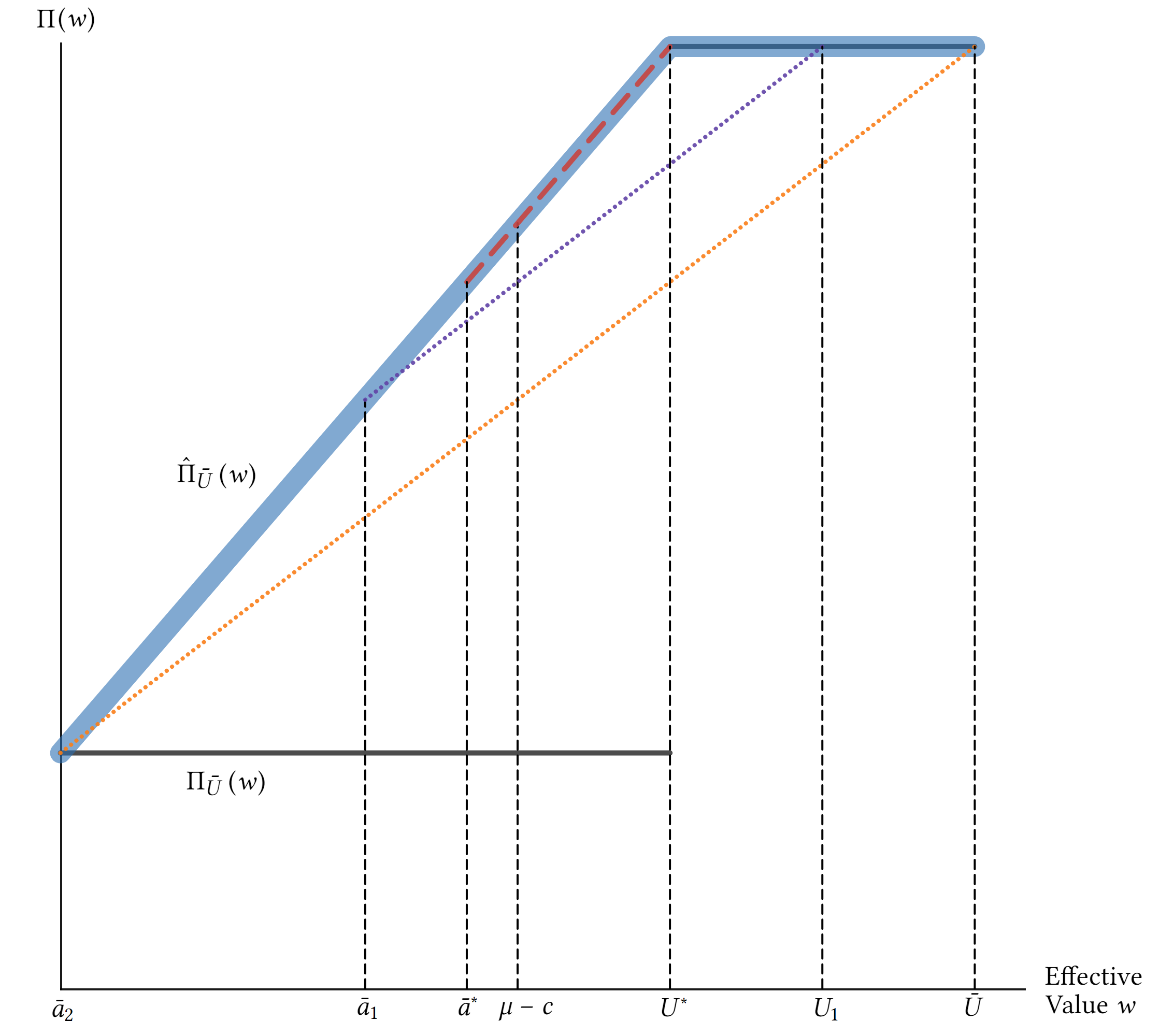}
    \caption{Example \ref{opt eff dist} analysis (for $\mu = 1/2$, $c = 1/8$, and $p_2 = 2/3$). The firm's payoffs from three different reservation values, $\bar{U}$, $U_1 = 5/8$, and $U^{*}$, are the intersections of $x = \mu - c$ and the dotted orange line, the dotted purple line, and the dashed red line, respectively. $U^{*}$ is clearly optimal.}
    \label{exfig}
\end{figure}

\begin{example}
\label{opt eff dist}Let's revisit the problem faced by firm 1 in Section
\ref{311}. Its payoff as a function of its realized effective value, $w$, is
\[\Pi \left(w\right) =
\begin{cases}
1-\frac{\mu }{p_{2}}, \quad &\text{if} \quad  w < p_{2}\bar{U} \\ 
1,\quad &\text{if} \quad w \geq p_{2}\bar{U}
\end{cases} \text{ .}
\]
It is clear that any effective-value distribution with reservation value falling short of $p_{2}\bar{U}$ results in a payoff of $1-\mu /p_{2}$. Take a reservation value $U\geq p_{2}\bar{U}$. The concavification of the restriction of the payoff function to domain $\left[ 0,U\right] $ is
\[\hat{\Pi}_{U}\left( w\right) =
\begin{cases}
1-\frac{\mu }{p_{2}}+\frac{\mu }{p_{2}^{2}\bar{U}}w, \quad &\text{if} \quad  w < p_{2}\bar{U} \\ 
1,\quad &\text{if} \quad w \geq p_{2}\bar{U}
\end{cases} \text{ .}
\]
From Proposition \ref{optimal prize distribution}, we need merely choose $U$ to maximize the following function.
\[V\left(U\right) = \left( 1-\frac{c}{1-U}\right) \times
\underbrace{\left( 1-\frac{\mu }{p_{2}}+\frac{\mu}{p_{2}\bar{U}}a\left(U\right)\right)}_{\hat{\Pi}_U\left(a\left(U\right)\right)} + \frac{c}{1-U} \times 1\text{ .}\]
Evidently, $V$ is strictly decreasing in $U$, so the optimal reservation value is $p_{2}\bar{U}$, and the corresponding optimal effective-value distribution has support $\left\{ 0,p_{2}\bar{U}\right\}$. Figure \ref{exfig} illustrates this analysis.
\end{example}

\section{The Symmetric Equilibrium in Competition over Effective Values}\label{eqsy}

In this section, we analyze the game of effective-value competition set up in the previous
section and characterize the unique symmetric equilibrium. We show that this equilibrium implies that each firm faces a payoff function that is linear in its effective value. We establish equilibrium existence by explicitly constructing a mixed strategy
that randomizes over binary distributions of posteriors.

Using the discrete choice reformulation outlined in Section \ref{secreform}, the interaction between firms can be modelled as competition over effective values. Formally, each of the $n$ firms simultaneously chooses an inducible effective-value distribution, with the objective of maximizing the probability that its realized effective value is the highest among those of all the firms (breaking ties fairly).

\begin{figure}
\centering
\begin{subfigure}{.5\textwidth}
  \centering
  \includegraphics[scale=.14]{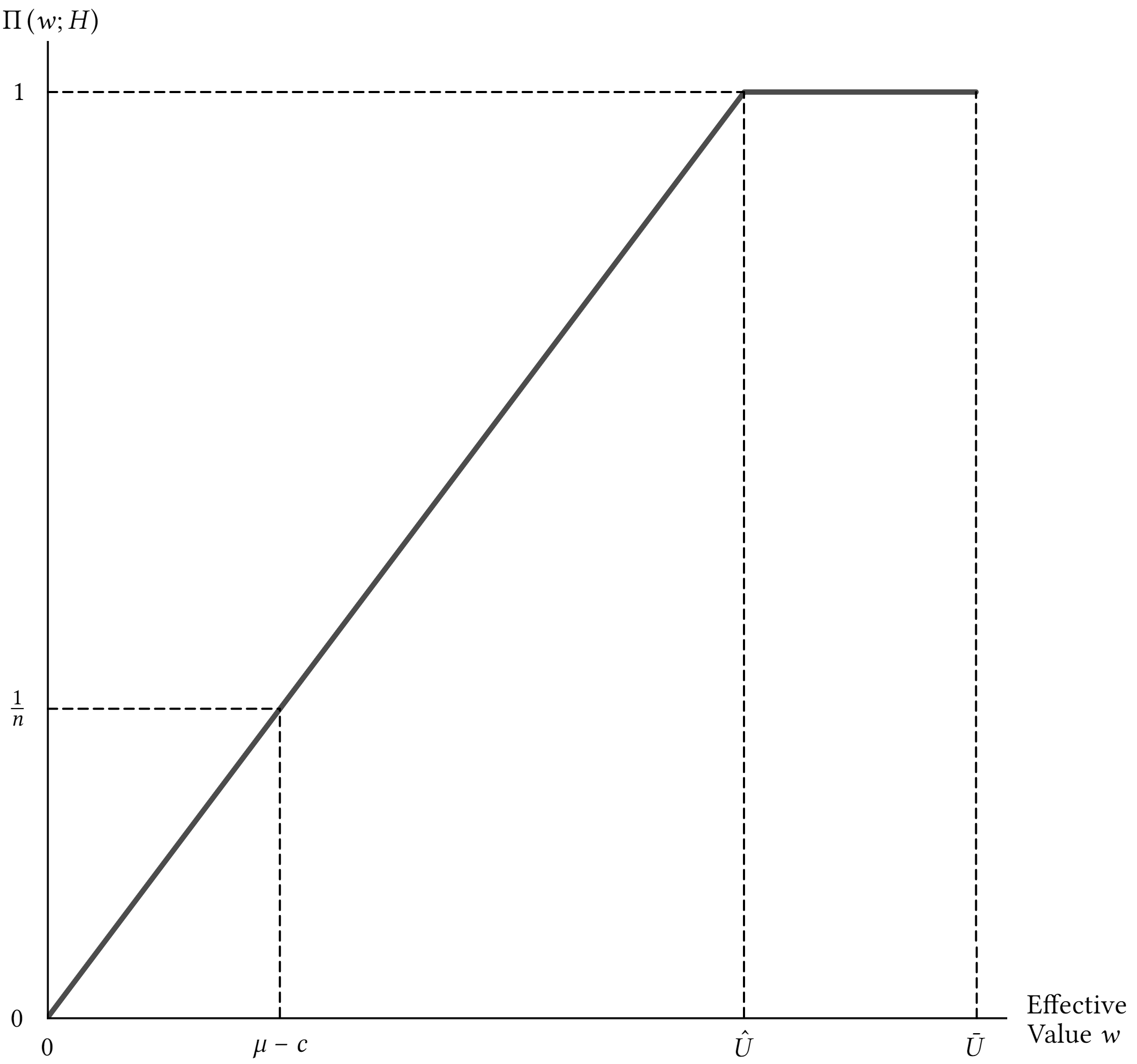}
  \caption{The equilibrium payoff when $\mu \leq \ubar{\mu}$.}
  \label{figsub1}
\end{subfigure}%
\begin{subfigure}{.5\textwidth}
  \centering
  \includegraphics[scale=.14]{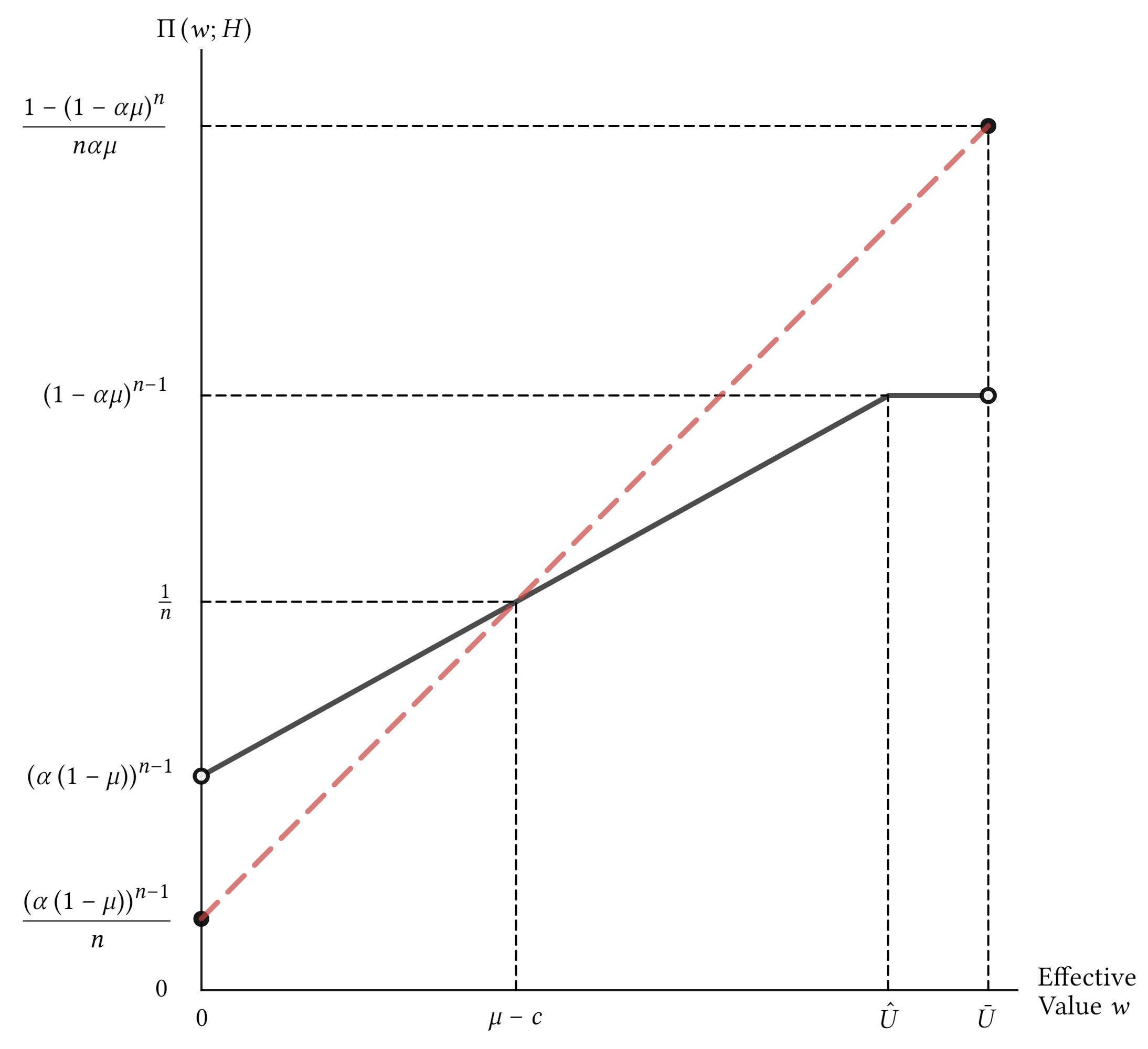}
  \caption{The equilibrium payoff when $\mu \in \left( \ubar{\mu},\bar{\mu}\right) $.}
  \label{figsub2}
\end{subfigure}
\par
\bigskip
\par
\bigskip
\par
\begin{subfigure}{.5\textwidth}
  \centering
  \includegraphics[scale=.14]{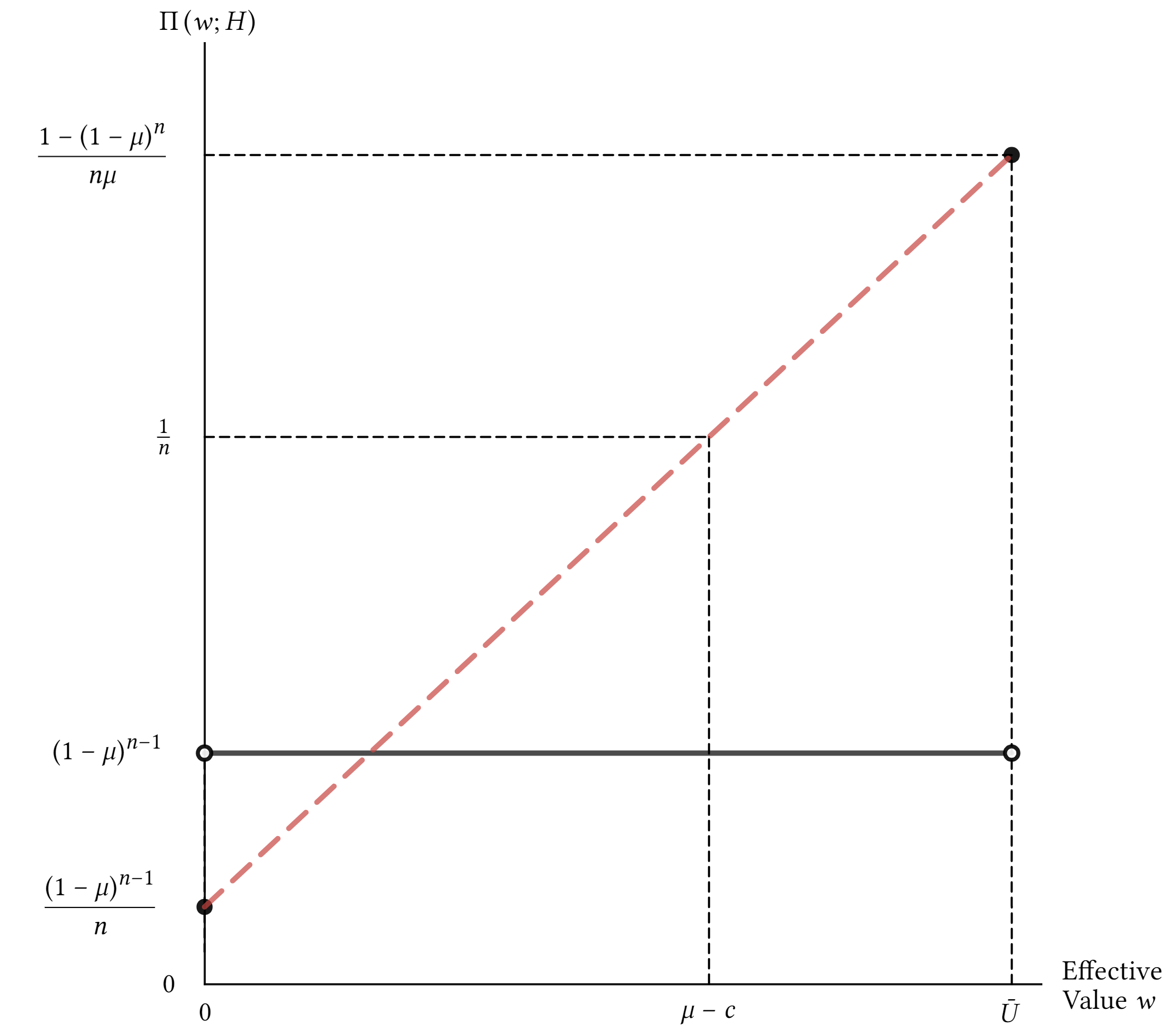}
  \caption{The equilibrium payoff when $\mu \geq \bar{\mu}$.}
  \label{figsub3}
\end{subfigure}
\caption{The linear structure of a firm's payoff function, $\Pi \left(
w;H\right) $.}
\label{fig2}
\end{figure}

We begin by noting that the equilibrium distribution $H_{i}$ must be atomless over $\left( 0,\bar{U}\right)$. To see this, observe that if the other firms place an atom at some $w\in \left( 0,\bar{U}\right)$, either as a reservation value or as a posterior realization, it is never optimal for a firm to respond by placing an atom there. If $w$ is a reservation value, offering a marginally more informative signal discretely improves the power of attraction. If $w$ is a posterior realization, shifting the weight to a marginally better posterior realization discretely improves the power of persuasion.\footnote{This argument does not apply to the effective value $0$ as full disclosure requires a mass there, a possibility that can arise in equilibrium (as we have shown in Section \ref{fullinfoeq} and will rederive below).} The next step in our equilibrium characterization is to establish that the payoff function $\Pi \left( w;H_{i}\right) $ facing each firm must possess the linear structure depicted in Figure \ref{fig2}. 

That the linear structure supports an equilibrium is straightforward. Indeed, the difficulty in showing that the linear structure corresponds to an equilibrium is in showing that it can be induced; which we do, constructively, in Lemma \ref{binary dist no u0}. In short, when an individual firm encounters a payoff function of the sort depicted in Figure \ref{fig2}, it is indifferent between offering any reservation value on the support of the payoff function and is thus willing to randomize over the relevant range. 

The linear structure is also necessary because the payoff function cannot have any strict convexities or concavities. Should we na\"{i}vely ignore the inducibility requirement for possible deviations, necessity would also be easy--straightforward concavification reveals that firms strictly prefer not to assign any positive mass to the interior of a region on which the payoff function is strictly convex and are unwilling to ``spread'' mass on any region on which the payoff function is strictly concave. That deviations must be inducible is relatively inconsequential when it comes to eliminating strict convexities. Alas, inducibility has more bite when showing that the equilibrium payoff cannot be strictly concave, and so we turn to the approach of Section \ref{modapproach}. Recall that a firm's best response can be characterized by its choice of reservation value $U$ and an optimal allocation of its flexible persuasion budget $a\left(U\right)$ (according to the concavification of $\Pi$ over $\left(0,U\right]$). Unless $\Pi$ is completely linear, the optimal allocation of $a\left(U\right)$ is constrained; e.g., if $\Pi$ is strictly concave on the interval in which $a\left(U\right)$ resides, it is not spread at all; and if it is locally linear, $a\left(U\right)$ can be spread only over those local values. We show that these implications inevitably lead to a contradiction; hence, $\Pi$ must be completely linear. 

We define $\ubar{\mu}\equiv 1/n$ and remind ourselves that $\Bar{\mu}\equiv 1-\left( 1/n\right)^{\frac{1}{n-1}}$. Then,
\begin{lemma}
\label{linear structure no u0}If $H_{i}$ is a symmetric equilibrium effective-value distribution, the implied payoff function in effective values facing each individual firm must have the following linear structure:
\[\Pi \left( w;H_{i}\right) = \begin{cases}
\frac{1}{n}\left( \alpha \left( 1-\mu \right) \right)^{n-1}, \quad &\text{if} \quad  w = 0 \\ 
\left( \alpha \left( 1-\mu \right) \right)^{n-1}+\frac{\left( 1-\alpha \mu
\right)^{n-1}-\left( \alpha \left( 1-\mu \right) \right)^{n-1}}{\hat{U}}w, \quad &\text{if} \quad w\in \left( 0,\hat{U}\right] \\
\left( 1-\alpha \mu \right)^{n-1}, \quad &\text{if} \quad w\in (\hat{U},\bar{U})\\
\frac{1-\left( 1-\alpha \mu \right)^{n}}{n\alpha \mu }, \quad &\text{if} \quad w=\bar{U}
\end{cases} 
\text{ ,}  \label{linear payoff no u0}\tag{$5$}\]
where $\alpha \in \left[ 0,1\right] $ and $\hat{U}\in \left[0,\bar{U}\right) $ are uniquely determined by $\mu $ and $n$ as follows.\footnote{This is this true except for the knife-edge case of $\mu =1/2$ and $n=2$. There, $\Bar{\mu}=\ubar{\mu}$ and cases \ref{casea}-\ref{casec} coincide. As a result, $\alpha $ can take any value in $\left[ 0,1\right] $ with $\hat{U}$ given in (\ref{alpha U_hat no u0 middle mean 2}).}
\begin{enumerate}[label={(\roman*)},noitemsep,topsep=0pt]
    \item\label{casea} If $\mu \leq \ubar{\mu}$, then $\alpha =0$ and $\hat{U}=n\ubar{U}$.
    \item\label{caseb} If $\mu \in \left(\ubar{\mu},\bar{\mu}\right)$, then $\alpha \in \left( 0,1\right) $ is the unique solution to \[\left( 1-\alpha \mu \right)^{n}-\left( \alpha \left( 1-\mu \right) \right)^{n}=1-\alpha \text{,}  \label{alpha U_hat no u0 middle mean 1}\tag{$6$}\]
    and \[\hat{U}=\frac{\left( 1-\alpha \mu \right)^{n-1}-\left( \alpha \left( 1-\mu\right) \right)^{n-1}}{n^{-1}-\left( \alpha \left( 1-\mu \right)\right)^{n-1}}\ubar{U}\text{.}  \label{alpha U_hat no u0 middle mean 2}\tag{$7$}\]
    \item\label{casec} If $\mu \geq \bar{\mu}$, then $\alpha =1$ and $\hat{U}=0$.
\end{enumerate}
\end{lemma}
In the lemma above, $\alpha $ is the probability that a firm offers full information, and $\hat{U}$ is the maximum reservation value, other than $\bar{U}$, that is on the support of a firm's equilibrium strategy. We will argue in Lemma \ref{binary dist no u0} that an equilibrium that generates the payoff function (\ref{linear payoff no u0}) exists. Before doing so, let us discuss the economic content of Lemma \ref{linear structure no u0}.

First, the linearity of the payoff function (\ref{linear payoff no u0}) highlights how the trade-off between attraction and persuasion in a firm's signal design problem is balanced in equilibrium. Choosing a more informative signal, and hence a high reservation value, $U$, caters to the attraction incentive but comes at the expense of persuasion effectiveness. The linearity of the payoff function (\ref{linear payoff no u0}) ensures that these forces cancel out, leaving the firm exactly indifferent between offering a range $\left[\ubar{U},\hat{U}\right] $ (and possibly $\bar{U}$) of reservation values. Moreover, not only is the firm indifferent as to what reservation value it induces, but it is also indifferent about how it allocates the flexible persuasion budget permitted by the chosen $U$. The rationale for this nonchalance is the usual matching pennies logic familiar from the papers on competitive persuasion without frictions (or papers with price dispersion): the firms randomize to leave the others indifferent and hence willing, themselves, to randomize. 

Second, the lemma shows that the equilibrium level of signal informativeness varies with the average product quality $\mu$ and the number of firms $n$. Figure \ref{fig2} depicts the three possible forms a firm's payoff function may take and thereby indicates the three possible guises of the equilibrium distributions. If the average quality and the number of firms are sufficiently low ($\mu \leq \ubar{\mu}$), firms provide partial information with probability one. On the other hand, if the average quality and the number of firms are sufficiently high ($\mu \geq \bar{\mu}$), firms provide full information with certainty. If the average quality and the number of firms are in some intermediate region ($\mu \in \left(\ubar{\mu},\bar{\mu}\right)$), both partial and full information coexist in equilibrium--it is a mixture of the other two cases.\footnote{Notably, ties occur with strictly positive probability (at equilibrium) outside of the low region; \textit{viz.}, if and only if $\mu > \ubar{\mu}$.}

As we discussed above, with few competitors (small $n$) and a low chance of the consumer stopping at
a competitor (small $\mu $), firms do not find it very difficult or important to attract the consumer. In particular, full information is overkill: even if a firm can capture the consumer's first visit, it can be converted into a sale only if its posterior realization is $1$, an event with low probability when $\mu$ is low. As a result, in equilibrium, no firms provide full information, and partial information (yielding a reservation value $\hat{U}=n\ubar{U}$) is sufficient to guarantee the top spot in the consumer's search order.

An increase in $\mu$ and $n$ corresponds to an increase in competition, which strengthens the attraction motive. The reaction of firms to this change is to provide more information, which means that firms have to provide even more information to be visited early. The equilibrium balance can be struck only if all firms ``spread out'' their effective-value distributions. As a result, the maximum reservation value increases in both $\mu $ and $n$. When $\mu $ and $n$ become high enough ($\mu = \ubar{\mu}$), full information is required to guarantee the consumer's first visit.

At even higher levels of $\mu $ and $n$, $\mu \in \left( \ubar{\mu},\bar{\mu}\right) $, the attraction consideration becomes even more important. Now, firms must provide full information with a strictly positive probability in equilibrium; otherwise, full information would strictly dominate any form of partial information. Moreover, as $\mu$ and $n$ increase, the atom at full information has to grow in order to balance the payoffs of full and partial disclosure. Eventually, as $\mu $ and $n$ become sufficiently large, $\mu \geq  \bar{\mu}$, the atom at full information becomes one: every firm provides full information. This is precisely the scenario discussed in Section \ref{fullinfoeq}.

With the payoff function uniquely pinned down by Lemma \ref{linear structure no u0}, it is straightforward to back out the candidate equilibrium effective-value distribution $H_{i}$ (Claim \ref{eqdistexplicit} in the appendix). The final step in our equilibrium construction is to show that the effective-value distribution thus obtained is indeed inducible; i.e., that it can be implemented via some mixed strategy over signals.

\begin{lemma}
\label{binary dist no u0}The effective-value distribution $H_{i}$
characterized in Lemma \ref{linear structure no u0} is inducible. Moreover,
it can be implemented by a mixed strategy that randomizes over binary
distributions of posteriors.
\end{lemma}

The proof of Lemma \ref{binary dist no u0} constructs a simple (not necessarily unique) mixed strategy that yields $H_{i}$. The binary implementation described in this lemma has several interesting properties. Once the initial randomness from the firms mixing is resolved, the consumer is faced with $n$ firms, each with binary distributions that nest within each other like a matryoshka doll. The experiments chosen by the firms can be ranked according to the Blackwell order, and the consumer searches them in order of their informativeness. The consumer stops only if she observes the high realization at a firm. Otherwise, she continues her search, and selects the last firm no matter its realization. Though this is a search in which recall is allowed, the consumer never utilizes this, and never returns to a firm from which she had previously moved on. 


This section's analysis is summarized by the following theorem.

\begin{theorem}
\label{summary no u0}Except for the knife-edge case in which $\mu =1/2$ and $%
n=2$, there exists a unique symmetric equilibrium in the firms' competition
over effective values. If $\mu =1/2$ and $n=2$, there exists a continuum of
symmetric equilibria. The equilibrium effective-value distribution gives
rise to a payoff function with a linear structure as in (\ref{linear payoff no u0}).
\end{theorem}

We conclude this section by discussing the effect of the search cost, $c$, on the consumer's welfare. While an increase in $c$ has a direct negative effect on the consumer's payoff, it induces more intense competition between the firms, which could potentially benefit the consumer. In fact, in the price-competition setting of \cite{choi2}, the indirect effect of intensified competition can be so strong that the consumer surplus increases with a higher search cost. In contrast, we find that when the search cost, $c$, is positive, the indirect effect only partially offsets the direct
negative effect, to the consumer's detriment.

\begin{corollary}
\label{CS_welfare}For $c>0$, an increase in $c$ worsens the equilibrium
distributions over effective values in the sense of first-order stochastic
dominance, thus hurting the consumer's \textit{ex ante} welfare.
\end{corollary}

\section{Two Benchmarks}

\label{bench}

In our main model, the consumer can learn the firms' posted signals for free, but discovering their realizations is costly. Consequently, a firm's posted signal plays the dual role of attracting the consumer to visit it and persuading her to purchase from it. The role of attraction is key in driving both our full disclosure and information dispersion results; and to illustrate this, we consider two alternative scenarios. In the first, the consumer can learn both the posted signals and their realizations for free. In the second, the signals are no longer posted, so the consumer only discovers a firm's signal and its realization \textit{after} paying the search cost. In contrast to our main model, neither setting generates information dispersion as the unique equilibrium outcome. Moreover, we show that equilibrium disclosure lessens, possibly radically, when the signals do not play a role in enticing the consumer to visit.

\subsection{Costless Signal Realizations}

\label{costlessbench}

This subsection considers the case without search frictions ($c=0$) and compares its equilibrium outcomes to those from the case with search frictions ($c>0$) studied in Section \ref{eqsy}. Without search frictions, the consumer has simultaneous and free access to the signal (posterior) realizations of each firm. While this frictionless case admits a host of symmetric equilibria depending on how ties are resolved, we continue to assume a fair tie-breaking rule; i.e., that the consumer randomizes equally between firms that offer the maximum posterior realizations. The justification of this equilibrium selection is as follows.\footnote{On top of the two reasons that we mention, it can be shown that equilibria that rely on the consumer's differential treatment of firms based on payoff-irrelevant histories (such as signal informativeness) are not robust to perturbations of the consumer's valuations for the sellers' products.} First, it enables a meaningful comparison. Recall that the equilibrium identified in Proposition \ref{summary no u0} also has the consumer adopt the fair tie-breaking rule, thereby following a Markov strategy that treats all firms equally (with the state variables being the firms' effective-value realizations). In the frictionless case, the fair tie-breaking rule remains a Markov strategy that treats all firms equally (with the state variables being the firms' posterior realizations).\footnote{We explain below why posterior realizations are natural counterparts of effective-value realizations in the frictionless case.} Second, even if the consumer's sampling of the signal realizations in the frictionless case is sequential, it is weakly dominant for her to sample all of them, as the visits are costless. This equilibrium serves as a useful benchmark because it eliminates the attraction motive altogether (since a firm's posterior offering will be observed irrespective of its signal choice).

Regardless of how small $c$ is, as long as it is positive, the consumer must discover the firms' signal realizations in sequence, which engenders the attraction motive. In contrast, when $c=0$, the consumer can costlessly collect and observe all signal realizations, and the attraction force vanishes completely. The qualitative difference between the two scenarios leads to drastically different equilibrium outcomes. First, with $c=0$, a symmetric pure-strategy equilibrium is possible. Second, the (symmetric) equilibrium distribution of effective values is not always continuous at $c=0 $. When it is not ($\mu >\ubar{\mu}$), the informativeness of the equilibrium signal jumps down discretely at $c=0$. This implies, somewhat counterintuitively, that the consumer may strictly prefer a small positive search cost to no search cost.

When $c=0$, equation (\ref{reservation value}) implies that all feasible distributions over posteriors have reservation values equal to one and hence, each implied effective value coincides with the posterior. As a result, the competition over effective values reduces to a straightforward competition over posterior realizations. This game is studied by \cite{Au2},\footnote{See also Section \ref{relwork}, which references other papers that investigate variants of the frictionless problem.} who show that a unique pure-strategy symmetric equilibrium exists. Intuitively, a pure-strategy equilibrium is possible here because overbidding one's rivals by offering a marginally more informative signal no longer results in an increase in the reservation value. A firm's signal, therefore, has no impact on the likelihood of being inspected by the consumer. As long as the payoff function in posteriors is linear, the firm is happy to put positive weight on the range of posteriors over which the payoff function is increasing. As a result, information dispersion is not a necessary feature
of the equilibrium.

The difference in nature between the limit game (with $c=0$) and the limiting game (with an arbitrarily small but positive $c$) raises the natural question of continuity, which is addressed in the proposition below.

\begin{proposition}
\label{continuity}Let $H_{i}^{c}$ be the equilibrium distribution of effective values when the search cost is $c\in \left[ 0,\mu \right] $, and let $H_{i}^{\ast }$ be the limiting distribution as the search cost vanishes, i.e., $H_{i}^{c}\rightarrow H_{i}^{\ast }$ in distribution as $c\rightarrow 0$.
\begin{enumerate}[label={(\roman*)},noitemsep,topsep=0pt]
    \item If $\mu \leq \ubar{\mu}$, $H_{i}^{\ast }=H_{i}^{0}$.
    \item If $\mu >\ubar{\mu}$, $H_{i}^{\ast }$ is a mean-preserving spread of $H_{i}^{0}$.
\end{enumerate}
\end{proposition}

It is immediate that the consumer's welfare increases if every firm provides more information. In addition, if $\mu > \ubar{\mu }$ the consumer may benefit from a small positive search cost $c$, which allows her to commit to reward informativeness by visiting firms with higher reservation values first. This commitment power intensifies the competition in disclosure between the firms, to the benefit of the consumer. In fact, if $\mu \geq \bar{\mu}$, then a positive search cost ensures full disclosure by the firms. In contrast to the classical Diamond paradox (\cite{diamond}); here, a small search cost begets the perfect competition (first-best) level of information provision.%

\subsection{Hidden Signals}

\label{diamond}

This subsection considers the case in which the firms' signals are not
directly observable to the consumer at the outset of her search. In
particular, she must incur search cost $c>0$ to discover both a firm's
signal and its realization. Accordingly, a firm's choice of signal cannot
affect the consumer's search order. Similar to the benchmark in the previous
subsection, this scenario is one of pure persuasion; and with the signal's role of attraction evaporated, it is natural to expect that the informativeness of firms' signals should decrease. In fact, we find that the
unobservability of signals poses a severe holdup problem akin to the Diamond paradox: each firm has an incentive to secretly lower the information content of its signal to increase the chance of successful persuasion once
the consumer has paid it a visit. This begets a stark equilibrium outcome: each firm's signal is uninformative, and the consumer visits only one firm.\footnote{The equilibrium concept appropriate for this benchmark is Perfect Bayesian Equilibrium, in which the consumer's conjectures about other theretofore unvisited firms' behavior are unaffected by observed deviations (``no signaling what you don't know''). We omit a formal definition of this for the sake of brevity.}

The reasoning is as follows. As no randomness is resolved until the consumer visits, it is without loss to focus on pure strategies of the firms. Consider a purported equilibrium in which some firms are believed to provide useful information to the visiting consumer, i.e., conditional on the visit, the probability that the consumer stops the search and purchases from the
visited firm is less than one. Let $\tilde{F}_{i}$ denote the consumer's conjecture of firm $i$'s distribution over posteriors and let $\tilde{U}_{i}$
be the corresponding reservation value. Given the consumer's (correct) belief about the equilibrium strategies of the firms, there exists a cutoff
posterior realization, denoted by $z_{i}^{\ast }$, above which the consumer stops the search and buys from firm $i$ conditional on visiting it.\footnote{If the purported equilibrium is symmetric, the value of this cutoff is common for all firms and is simply $\tilde{U}_{i}$. If the purported equilibrium is asymmetric, this cutoff is the optimal continuation value of search beyond firm $i$, which is weakly below $\tilde{U}_{i}$.} It follows from \hyperlink{pandora}{Pandora's rule} that $z_{i}^{\ast }\leq 
\tilde{U}_{i}$. If $z_{i}^{\ast }\leq \mu $, then firm $i$ can secure the
consumer's purchase (conditional on visit) by a distribution over posteriors
supported on $\left[ z_{i}^{\ast },1\right] $; that is, firm $i$ provides no
useful information to the consumer. Conversely, if $z_{i}^{\ast }\in (\mu ,%
\tilde{U}_{i}]$, then the optimal signal of firm $i$ assigns no weight to
posteriors above $z_{i}^{\ast }$,\footnote{%
This is a simple consequence of the concavification technique of \cite{kam}.}
making reservation value $\tilde{U}_{i}$ impossible. As a result, the only equilibria involve all firms providing useless information to the consumer, who stops her search at the first firm. The following theorem summarizes the discussion above.

\begin{theorem}
\label{infodia} Suppose the firms' choice of signals are revealed to the consumer only after she pays the visit cost. In all equilibria, each firm
offers a distribution of effective values that is degenerate at $\ubar{U}$, and the consumer buys from the first visited firm with probability one. 
\end{theorem}

The result can be understood as an informational Diamond paradox--in all equilibria, firms provide only the monopoly level of information, and there is no active consumer search. Without the attraction motive, the persuasion motive encourages pooling above the consumer's stopping threshold, which yields a profitable deviation to any firm that is providing useful information. That the inability of firms to secretly deviate and provide less information than anticipated is important for sustaining equilibria with (useful) information provision is a robust finding. \cite{board} establish a similar ``informational Diamond paradox'' in their Theorem 2; and \cite{dia} shows that this insight persists even with price setting.  

\section{Extensions}

\label{extensions}

In this section, we discuss several extensions of the basic model. The details and formal results are left to the \href{https://whitmeyerhome.files.wordpress.com/2022/05/attraction_versus_persuasion_supplementary_appendix_jpe_final.pdf}{Supplementary Appendix}.

First, the \textit{ex ante} homogeneity of the firms justifies our focus on symmetric equilibria. What about asymmetric equilibria? There are none when there are just two firms, but for three or more firms there is one when $\mu$ is smaller but sufficiently close to $\bar{\mu}$. In this equilibrium, $n-1$ firms provide full information, and the $n$th firm chooses a distribution over effective values that is identical to the equilibrium distribution when there are two firms and $\mu$ is low. This equilibrium is notable as it yields higher consumer welfare than the coexisting symmetric equilibrium--in fact, firms provide the first-best amount of information.

Second, we also explore the scenario with two heterogeneous firms (\textit{viz.}, with different expected qualities). There are four different regions of the parameter space, each of which begets a different variety of equilibrium. If the gap between the means is large enough, the higher mean firm is visited first and selected in all equilibria. In the other regions, the familiar linear structure manifests, although one or both firms may place masses on either $0$ or the maximal reservation value the lower mean firm can induce.

Third, we allow the consumer a positive outside to which she may always return upon quitting her search. Mirroring Lemma \ref{linear structure no u0}, there is a symmetric equilibrium that gives rise to a linear payoff function. The relevant outside option means that the game is no longer zero sum: surprisingly, industry profits can be hurt not only by an increase in the consumer's search cost $c$,\footnote{\cite{choi2} report a similar finding in their price-competition setting.} but also by an improvement in the average product quality $\mu$. A higher average quality incites more aggressive information revelation (as attraction becomes more important), which harms market profits by increasing the probability that the consumer takes her outside option.

Fourth, we illustrate that the important qualitative features of the earlier results do not rely on the binary prior and, in particular, persist when the consumer's match value is distributed according to some absolutely continuous prior. The tension between attraction and persuasion remains: if the number of firms in the market is sufficiently high, the attractive motive is strong, and so any pure strategy equilibrium must be one in which firms compete fiercely to be visited first and, therefore, induce the maximum reservation value. Otherwise, no symmetric equilibrium exists; instead, firms randomize over signals.

\section{Discussion and Concluding Remarks}

\label{discuss}

\label{dis}

This paper explores competition in information provision in a sequential, directed search setting. Our model highlights the key economic forces at work in this environment; namely, the conflicting motives of attraction and persuasion. We uncover a number of insights pertaining to how the underlying environment shapes these two incentives and the ensuing information provision. For instance, a sufficiently high average quality or a large number of competitors makes the attraction motive dominant, leading to full
information in equilibrium. 

Outside of the high average quality case, although the attraction incentive remains, the persuasion motive has more of an effect: if everyone else provides full information, it is now worthwhile for a firm to provide no information and count on the consumer to visit and select it at the end of an unsuccessful search. The forces of attraction and persuasion can be balanced only in a mixed-strategy equilibrium, resulting in information dispersion. By contrasting this finding with alternative settings in which the attraction motive is irrelevant, we show not only that the attraction motive is key to generating the information dispersion result, but also, somewhat counterintuitively, that the consumer can actually benefit from having a small positive search cost (rather than none).

The assumption that firms have perfect flexibility in their ability to design signals not only aids us in obtaining a straightforward equilibrium characterization but also helps make transparent the aforementioned tension between attraction and persuasion. While reducing noise always helps attract the consumer, the effectiveness of persuasion can be enhanced only by
introducing noise in a specific manner. A restrictive set of feasible signals might thus obscure our model's essential trade-off. Moreover, although in practice, a firm's control over information revelation is never
perfect, the economic insights we uncover remain valid provided the set of feasible signals is not too meagre. Consider, for instance, the result of full disclosure with a sufficiently high $\mu $ or $n$. If the signal space were more restrictive than the one we consider, the attraction motive would still dominate if competition were sufficiently intense. 

Our analysis provides a number of testable predictions that await empirical investigation. One of the main ideas emerging from our study is the possibility of information dispersion when competition is not that intense. Do we observe such information dispersion in real-world markets? Unlike price dispersion, which is relatively easy to observe and measure (since prices are scalars), information levels are much harder to quantify, which could explain the dearth of formal evidence of this phenomenon. This difficulty notwithstanding, casual observation suggests that such variation does exist. An alternative interpretation of our model is that firms compete by choosing their product designs, which affect the distributions of their match values with the consumers. Naturally, a broad design induces a more concentrated distribution of match values, whereas a niche design induces a more spread-out distribution. With this interpretation, our model suggests that design dispersion can arise despite  \textit{ex ante} homogeneity. 

\newpage

\appendix

\section{Sections \protect\ref{prelim}, \protect\ref{eqsy}, and \protect\ref%
{bench} Proofs}

\subsection{Proof of Lemma \protect\ref{no partial pure eqm}}

This is a special case of Lemma \ref{no mass atom} and is thus omitted.

\subsection{Proof of Proposition \protect\ref{2states}}

This proposition is a special case of Theorem \ref{summary no u0};
specifically, case \ref{casec} of Lemma \ref{linear structure no u0}.

\subsection{Proof of Lemma \protect\ref{effective-value mean}}

\textit{Proof.}
\begin{enumerate}[label={(\roman*)},noitemsep,topsep=0pt]
    \item The existence of an atom at $U$ and that $H\left( U;F_{i}\right) =1$ follow from (\ref{eff-value dist.}). It remains to show that the atom size is at least $c/\left( 1-U\right) $. By (\ref{reservation value}), $c=\left( 1-F\left( U\right) \right) \times \left[ \int_{U}^{1}\left( x-U\right) dF_{i}\left( x\right) /\left(1-F\left( U\right) \right) \right] $. As the term in brackets is no larger than $1-U$, we have $1-F\left( U\right) \geq c/\left( 1-U\right) $.
    \item Given a distribution $F_{i}$ over posteriors, the mean of the implied effective-value distribution $H\left( \cdot ;F_{i}\right) $ is 
    \[\begin{split}
    \int_{0}^{\bar{U}}wdH\left(w; F_i\right) &= \bar{U}-\int_{0}^{\bar{U}}H\left(w; F_i\right) dw = \bar{U}-\int_{0}^{U}F_i\left( w\right) dw-\int_{U}^{\bar{U}}1dw \\
    &= U-\left( U+c-\mu \right) =\ubar{U}
    \end{split} \text{ .}\]
We obtain the first equality via integration by parts and the second equality through the definition of $H\left( \cdot ;F_{i}\right) $ (given by (\ref{eff-value dist.})). The third equality uses the reservation-value equation (\ref{reservation value}), which implies $\int_{0}^{U}F_{i}\left( w\right) dw=U+c-\mu $. \hfill $\blacksquare$
\end{enumerate}

\subsection{Proof of Lemma \protect\ref{inducibility}}
\begin{proof}
Take a mixed strategy $\sigma_{i}\in \Delta \left( \mathcal{F}\right) $ of
firm $i$. Let $G_{\sigma _{i}}\left( u\right) \equiv \sigma _{i}\left(\left\{ F_{i} \colon U\left( F_{i}\right) \leq u\right\} \right) $ be the implied distribution of reservation values. Let $\sigma _{i}\left( \cdot \vert u\right) \in \Delta \left( \mathcal{F}\right) $ be the conditional distribution over
pure strategies that $\sigma _{i}$ induces given reservation value $u\in %
\left[\ubar{U},\bar{U}\right]$.\footnote{The existence of this conditional distribution is a consequence of the Radon-Nikodym Theorem, see for instance, Section 33 of \cite{billingsley}.}
Furthermore, define, for each $u\in \left[ \ubar{U},\bar{U}\right] $, $%
F_{\sigma _{i}}\left( p;u\right) \equiv \int F_{i}\left( p\right) d\sigma
_{i}\left( F_{i} \vert u\right) $ the distribution over posteriors implied by $\sigma _{i}$ conditional on reservation value $u$. The effective-value
distribution $H_{i}\left(w;\sigma _{i}\right) $ induced by $\sigma _{i}$
can then be expressed as
\[\begin{split}
    H_{i}\left( w;\sigma _{i}\right) = \Pr \left( \min \left\{ p_{i},U\left(
F_{i}\right) \right\} \leq w\right) &= \int \left( \mathbf{1}_{\left[ U\left( F_{i}\right) \leq w\right]
}+F_{i}\left( w\right) \mathbf{1}_{\left[ U\left( F_{i}\right) >w\right]
}\right) d\sigma _{i}\left( F_{i}\right) \\
&= G_{\sigma _{i}}\left( w\right) +\int_{w}^{\bar{U}}\underset{F_{\sigma
_{i}}\left( w;u\right) }{\underbrace{\left[ \int F_{i}\left( w\right)
d\sigma _{i}\left( F_{i} \vert u\right) \right] }}dG_{\sigma _{i}}\left( u\right)
\end{split} \text{ .}\]

Conversely, given a reservation-value distribution $G_{i}\in \Delta \left( %
\left[ \ubar{U},\bar{U}\right] \right) $ and a feasible
distribution over posteriors $F_{i}\left( \cdot ;u\right) \in \mathcal{F}$
for each reservation value $u$ on the support of $G_{i}$, one can construct
the corresponding mixed strategy $\sigma \in \Delta \left( \mathcal{F}\right) $ by $\sigma \left( B\right) =\int_{B\cap \left\{ F_{i}\left( \cdot ;u\right) \colon u\in \supp\left\{ G_{i}\right\} \right\} }dG_{i}\left( u\right)$ for each $B\subset \mathcal{F}$ Borel-measurable. \end{proof}

\subsection{Proof of Lemma \protect\ref{linear structure no u0}}
\begin{proof}
Here is a roadmap for the proof:
\begin{enumerate}[noitemsep,topsep=0pt]
    \item Claim \ref{no mass atom} eliminates the possibility of atoms at effective values other than $0$ or $\bar{U}$. 
    \item Claim  \ref{full info unique} determines the value of $\alpha$, the probability that firms provide full information.
    \item Claim \ref{U hat} derives the value of $\hat{U}$ specified in the lemma's statement and verifies that it is feasible.
    \item Claim \ref{concavity} reveals the necessary (weak) concavity of the equilibrium payoff function on $\left[0,\ubar{U}\right]$.
    \item Claim \ref{strictconcavity} argues that the equilibrium payoff function cannot be strictly concave.
    \item Claim \ref{nomaximalint} establishes that the equilibrium payoff function cannot be piecewise linear. 
    \item Claim \ref{eqdistexplicit} unearths the equilibrium distribution over effective values.
\end{enumerate}

\begin{claim}[No Atoms Except Possibly at $0$ and/or $\bar{U}$]
\label{no mass atom} In any symmetric equilibrium, a firm's effective-value distribution $H_{i}$ has no atoms except possibly at $0$ and $\bar{U}$. Consequently, the expected payoff facing an individual firm $i$, as a function of its realized effective value $w$, takes the following form
\[\Pi \left( w;H_{i}\right) \equiv 
\begin{cases}
\frac{1}{n}H_{i}\left( 0\right)^{n-1}, \quad &\text{if} \quad  w = 0 \\ 
H_{i}\left( w\right)^{n-1}, \quad &\text{if} \quad w\in (0,\bar{U})\\
\lim_{w^{\prime }\rightarrow \bar{U}^{-}}\frac{1-H_{i}\left( w^{\prime
}\right)^{n}}{n\left( 1-H_{i}\left( w^{\prime }\right) \right) },\quad &\text{if} \quad w=\bar{U}
\end{cases} 
\text{ .}  \label{payoff in H}\tag{$A1$}\]
\end{claim}
\begin{proof}
Keep in mind that the effective-value distribution $H_i$ is derived using the discrete choice reformulation of \cite{am} and \cite{choi2}, the construction of which is explained in Lemma \ref{inducibility}. Suppose $H_i$ has an atom at some $\tilde{w}\notin \left\{ 0,\bar{U}\right\} $. We show that the best response to $\Pi \left( \cdot ;H_i\right) $ does not
put any positive mass at $\tilde{w}$. Proposition \ref{optimal prize
distribution} implies the best response to $\Pi \left( \cdot ;H_i\right) $ puts a
positive mass at $\tilde{w}$ only if either there is some $U>\tilde{w}$ such
that $\Pi\left( \tilde{w};H_i\right) =\hat{\Pi}_{U}\left(\tilde{w};H_i\right) $, or if $\tilde{w}$ is an optimal reservation value. First, an atom of $H_i$ at $\tilde{w}$ implies $\Pi \left( \tilde{w};H_i\right) <\lim_{w\rightarrow \tilde{w}^{+}}\Pi \left( w;H_i\right) $. Consequently, $\Pi\left(\tilde{w};H_i\right) <\hat{\Pi}_{U}\left(\tilde{w};H_i\right) $
for all $U>\tilde{w}$, so at every reservation value $U\in \left[\ubar{U},\bar{U}\right] $, it is suboptimal to assign positive mass at effective value $\tilde{w}$ (recall the concavification approach of Section \ref{modapproach}). As a result, $\Pi \left(w;H_{i}\right) $ must be continuous at all $w\in \left(0,\ubar{U}\right]$.

Next, $\tilde{w}$ cannot be an optimal reservation value either. To see this, recall the maximum payoff given reservation value $\tilde{w}$ is 
\[V\left( \tilde{w}\right) =\left( 1-\frac{c}{1-\tilde{w}}\right) \times \hat{%
\Pi}_{\tilde{w}}\left( a\left( \tilde{w}\right) ;H_{i}\right) +\frac{c}{1-%
\tilde{w}}\times \Pi \left( \tilde{w};H_{i}\right) \text{ .}\]
As $\tilde{w}\neq \bar{U}$, $a\left(\tilde{w}\right) \in \left(0,\ubar{U}\right]$, and we know from above that $\hat{\Pi}_{\tilde{w}}\left( \cdot ;H_{i}\right) $ is continuous at $a\left( \tilde{w}\right) $. This plus the continuity of $a\left( \cdot \right) $ implies $\hat{\Pi}_{w}\left( a\left( w\right) ;H_{i}\right) $ is continuous at $w=\tilde{w}$. The atom of $H_{i}$ at $\tilde{w}$ implies $\Pi \left( \tilde{w}; H_{i}\right) <\lim_{w\rightarrow \tilde{w}^{+}}\Pi \left(w;H_{i}\right)$; and consequently, $V\left( \tilde{w}\right) <\lim_{w\rightarrow \tilde{w}^{+}}V\left( w\right) $. \end{proof}
Recall that $\alpha \equiv \left(
1-\lim_{w\rightarrow \bar{U}^{-}}H_{i}\left( w\right) \right) /\mu $ is the probability firm $i$ offers full information.
\begin{claim}[Determining the Value of $\alpha$]
\label{full info unique}
\begin{enumerate}[label={(\roman*)},noitemsep,topsep=0pt]
    \item\label{firstunique} If $\mu \geq \bar{\mu}$, the symmetric equilibrium has $\alpha =1$.
    \item\label{secondunique} If $\mu \in \left(\ubar{\mu },\bar{\mu}\right) $, the symmetric equilibrium has $\alpha$ determined by (\ref{alpha U_hat no u0 middle mean 1}).
    \item\label{thirdunique} If $\mu \leq \ubar{\mu}$, the symmetric equilibrium has $\alpha =0$.
\end{enumerate}
\end{claim}
\begin{proof}
If each other firm chooses full information with probability $\alpha$, a firm's payoff from full information is
\[\mu \Pi \left(\bar{U};H_{i}\right) +\left( 1-\mu \right) \Pi\left(0;H_{i}\right) = \begin{cases}
\mu, \quad &\text{if} \quad  \alpha = 0 \\ 
\frac{1}{n}-\frac{1}{n\alpha }T\left( \alpha \right), \quad &\text{if} \quad \alpha \in \left(0,1\right)\\
\frac{1}{n},\quad &\text{if} \quad \alpha = 1
\end{cases} \text{ ,}\]
where $T\left( \alpha \right) \equiv \left( 1-\alpha \mu \right)
^{n}-\left( \alpha \left( 1-\mu \right) \right)^{n}-\left( 1-\alpha \right)$. The following observations will be useful: first, $0$ and $1$ are both roots of $T$. Second, $T^{\prime \prime }\left(\alpha \right) >0\Leftrightarrow \alpha <\left( \mu +\left( \mu^{-2}\left(1-\mu \right)^{n}\right)^{\frac{1}{n-2}}\right)^{-1}$, so $T^{\prime \prime }$ changes sign at most once, and the change can only be from
positive to negative. Third, $T^{\prime }\left( 0\right) =1-n\mu $ and $T^{\prime }\left( 1\right) =1-n\left( 1-\mu \right)^{n-1}$. In the special case in which $n = 2$ and $\mu = \bar{\mu} = \ubar{\mu} = 1/2$, a firm's payoff from full information is $1/2$ irrespective of the other firm's choice of $\alpha$. This yields the multiplicity described in the text. For the rest of the proof, assume $n \geq 3$.

\ref{firstunique} Let $\mu \geq \bar{\mu}$. If $\alpha =0$, the full-information payoff is $\mu$, which exceeds $1/n$ whenever $\mu \geq \bar{\mu}$. Because $\mu \geq \bar{\mu}$, $T^{\prime }\left( 0\right) <0$ and $T^{\prime }\left( 1\right) \geq 0$. Because $T^{\prime \prime }$ changes sign only once, $T^{\prime }$ can change sign only once (from negative to positive), and so $T\left( \alpha \right) <0$ for all $\alpha \in \left( 0,1\right)$. This yields a full-information payoff exceeding $1/n$, which is impossible. Thus, $\alpha$ must equal $1$.

\ref{secondunique} Let $\mu \in \left(\ubar{\mu },\bar{\mu}\right)$. Again, if $\alpha = 0$ the full-information payoff is $\mu > 1/n$, which is impossible. $\alpha =1$ at equilibrium is also infeasible, as it implies a full-information payoff of $1/n$, which is strictly below the no-information payoff of $\left( 1-\mu \right)^{n-1}$ whenever $\mu <\bar{\mu}$. This leaves $\alpha \in \left(0,1\right) $ as the only possibility. Because $\mu \in \left( \ubar{\mu},\bar{\mu}\right) $, $T^{\prime }\left( 0\right) <0$ and $T^{\prime }\left( 1\right) <0$. As $T\left( 0\right) =T\left( 1\right) =0$, the continuity of $T\left( \alpha \right) $ implies that there must exist an interior root. Moreover, as $T^{\prime \prime }\left( \alpha \right) $ changes sign only once, there exists at most one interval over which $T^{\prime }\left( \alpha \right) $ is positive. The interior root of $T\left(\alpha \right) $ is thus unique, and $T\left( \alpha
\right) =0$ is equivalent to (\ref{alpha U_hat no u0 middle mean 1}).

\ref{thirdunique} Let $\mu \leq \ubar{\mu}$. $\alpha =1$ is impossible at equilibrum, as it implies a full-information payoff of $1/n$, which is strictly below the no-information
payoff of $\left( 1-\mu \right)^{n-1}$ whenever $\mu \leq \ubar{\mu }< \bar{\mu}$. $\mu \leq \ubar{\mu}$ implies $T^{\prime }\left(0\right) \geq 0$ and $T^{\prime }\left( 1\right) <0$; and because $T^{\prime \prime }$ changes sign only once $T^{\prime }$ also changes sign only once. As $T\left( 0\right) =T\left( 1\right) =0$, $T\left( \alpha \right) >0$ for all $\alpha \in \left(0,1\right) $. Accordingly, an equilibrium cannot have a nonzero $\alpha $; otherwise, the full-information payoff is below $1/n$, a contradiction. This leaves $\alpha =0$ as the only possibility in equilibrium.
\end{proof}
\begin{claim}[Determining the Value of $\hat{U}$]
\label{U hat}In any symmetric equilibrium with a linear payoff function, the
value of $\hat{U}$ is $\left(\supp\left( H_{i}\right) /\left\{ \bar{U}\right\} \right) $, as specified in the lemma's statement. Moreover, $\hat{U}$ is feasible.
\end{claim}

\begin{proof}
The case of $\mu \geq \bar{\mu}$ is immediate, so we focus on the other two
cases.

If $\mu \leq \ubar{\mu}$, Claim \ref{full info unique} states that $\alpha =0$. Moreover, because the linear payoff function facing each firm implies a payoff of $1/n$ at equilibrium, we must have $\Pi\left(\ubar{U};H_{i}\right) =1/n$. This is possible only if $\hat{U}=n\ubar{U}$. Evidently, $\mu \leq \ubar{\mu}$ guarantees $\hat{U}\leq \bar{U}$.

Finally, if $\mu \in \left( \underaccent{\bar}{\mu },\bar{%
\mu}\right) $,  Claim \ref{full info unique} states that $\alpha \in \left( 0,1\right) $. Again, that the equilibrium payoff is $1/n$ requires $\Pi \left(\ubar{U};H_{i}\right) =1/n$ and 
\[\frac{\Pi \left( \ubar{U};H_{i}\right) -\lim_{w\rightarrow
0^{+}}\Pi \left( w;H_{i}\right) }{\ubar{U}}=\frac{\Pi \left( 
\hat{U};H_{i}\right) -\lim_{w\rightarrow 0^{+}}\Pi \left( w;H_{i}\right) }{%
\hat{U}}\text{ ,}\]
which upon rearranging gives (\ref{alpha U_hat no u0 middle mean 2}). It
remains to show that this value of $\hat{U}$ lies between $%
\ubar{U}$ and $\bar{U}$. $\hat{U}\geq \ubar{U}$ is equivalent to $n^{-1}\leq \left( 1-\alpha \mu\right) ^{n-1}$, which follows immediately from $\mu <\bar{\mu}$ and $\alpha \in \left( 0,1\right) $. Straightforward algebra reveals that $\hat{U}\leq \bar{U}$ is equivalent to $\left( 1-\alpha \mu \right) ^{n-1}\leq \alpha/n +\left( 1-\alpha \right) $. As $\alpha $ and $\mu $ are related by (\ref{alpha U_hat no u0 middle mean 1}), the last inequality is equivalent to requiring that for all $\alpha \in \left( 0,1\right) $,  
\[
\left( \alpha n^{-1}+\left( 1-\alpha \right) \right) ^{\frac{n}{n-1}}-\left(
1-\alpha \right) -\left( \left( \alpha n^{-1}+\left( 1-\alpha \right)
\right) ^{\frac{1}{n-1}}-\left( 1-\alpha \right) \right) ^{n}\geq 0\text{ .}
\]
To show this, denote $t\left( \alpha \right) \equiv \left( \alpha
n^{-1}+\left( 1-\alpha \right) \right) ^{\frac{1}{n-1}}$, so that $t$ has a
range of $\left( n^{-\frac{1}{n-1}},1\right) $. With this change of
variable, the inequality above can be expressed as 
\[S\left( t\right) \equiv \frac{\left( t-\frac{1}{n-1}\left( nt^{n-1}-1\right)
\right) ^{n}}{t^{n}-\frac{1}{n-1}\left( nt^{n-1}-1\right) }\leq 1 \quad \text{for
all} \quad t\in \left( n^{-\frac{1}{n-1}},1\right) \text{ .}\]
Evidently, the sign of $S^{\prime }\left( t\right) $ is equal to that of $-\left( n-2\right) t^{n-1}+\left( n-1\right) t^{n-2}-1$, which is clearly negative for all $t<1$, as it is increasing in $t$ and equal to $0$ at $t=1$. Therefore, $S\left( t\right) \leq S\left( n^{-\frac{1}{n-1}}\right) =1$.
\end{proof}

Next, we take the first step toward establishing the necessary linearity of the equilibrium payoff function, $\Pi \left(w;H_{i}\right) $, by showing that it must be weakly concave over the interval $\left[0,\ubar{U}\right]$.
\begin{claim}
\label{concavity}Suppose $\mu <\bar{\mu}$. Each firm can obtain the equilibrium payoff by offering no information. Moreover, $\Pi \left(w;H_{i}\right) $ is weakly concave over the interval $\left[0,\ubar{U}\right]$.
\end{claim}

\begin{proof}
To prove the first statement, it suffices to show $\Pi \left( \ubar{U}%
;H_{i}\right) =1/n$ in equilibrium. Suppose not. An atom must be found at $\inf \left(\supp\left( H_{i}\right) \cap \left(\ubar{U},\bar{U}\right]\right) $, contradicting the fact that $\Pi \left(w;H_{i}\right) $ has no interior atom.

The second statement is equivalent to $\Pi _{\ubar{U}}\left(
w;H_{i}\right) =\hat{\Pi}_{\ubar{U}}\left( w;H_{i}\right) $ for all $%
w\in \left(0, \ubar{U}\right]$.  Suppose to the contrary that there is some $w^{\prime }\in (0,\ubar{U})$ such that $\Pi\left( w^{\prime };H_{i}\right) <\hat{\Pi}_{\ubar{U}}\left( w^{\prime };H_{i}\right) $. Then there must be an open neighbourhood around $w^{\prime }$ over which no optimal effective-value distribution assigns any positive measure, and $\Pi\left(
\cdot ;H_{i}\right) $ must be flat over this neighbourhood. As $\Pi\left( w^{\prime };H_{i}\right) <\hat{\Pi}_{\ubar{U}}\left( w^{\prime };H_{i}\right) $ implies that there is some $w^{\prime
\prime }\in \left( w^{\prime },\ubar{U}\right]$ such that $\Pi\left( w^{\prime };H_{i}\right) <\Pi\left(w^{\prime \prime };H_{i}\right) $, $\Pi\left(\cdot;H_{i}\right)$ must have a discrete jump (and hence $H_{i}$ has an atom) somewhere in the interval $\left(w^{\prime },w^{\prime \prime }\right]$, contradicting Claim \ref{no mass atom} above.
\end{proof}

Note that an immediate implication of the weak concavity of $\Pi \left(
w;H_{i}\right) $ on $\left[0,\ubar{U}\right]$ is that value $0$ is
on the support of $H_{i}$. To establish the linearity of the payoff function
on $\left( 0,\hat{U}\right) $, it suffices to show that it is
linear on $\left( 0,\ubar{U}\right) $.\footnote{The reason is as follows. Suppose $\Pi \left( w;H_{i}\right) =L+\left(n^{-1}-L\right) w/\ubar{U}$ for all $w\in \left(0, \ubar{U}\right]$, where $%
L=\lim_{w\rightarrow 0}\Pi _{\ubar{U}}\left( w;H\right) $. If some $%
w^{\prime }\in \left( \ubar{U},\hat{U}\right) $ has $\Pi \left(
w;H_{i}\right) >L+\left( n^{-1}-L\right) w/\ubar{U}$, then some binary
effective-value distribution with $\left\{ a\left(w\right) ,w\right\} $ achieves a payoff exceeding $1/n$. If some $w^{\prime }\in \left( 
\ubar{U},\hat{U}\right) $ has $\Pi \left( w;H_{i}\right) <L+\left(
n^{-1}-L\right) w/\ubar{U}$, then any effective-value distribution with
this reservation value cannot yields a payoff less than $1/n$ and
thus cannot be on the support of $H_{i}$. Because $H_{i}$ has no
interior atom, it must be flat on $\left(w^{\prime },\bar{U}\right) $, contradicting the definition of $\hat{U}$.} Concavity implies that the left and right derivatives of $\Pi \left(w,H_{i}\right) $ exist at all points $(0,\ubar{U})$ and they can disagree at at most countably many points. Let $\bar{w}\equiv \sup\left\{\supp\left(H_{i}\right) \cap \left[0,\ubar{U}\right] \right\}$ (note $\bar{w}=\ubar{U}$ if $\ubar{U}\in \supp\left( H_{i}\right) $). If $\Pi\left( w,H_{i}\right)$ is concave, but not linear, on $(0,\ubar{U}]$, there are two possibilities:

\begin{enumerate}[noitemsep,topsep=0pt]
    \item $\Pi \left( w,H_{i}\right) $ is non-linear over $\left[ \bar{w}-\varepsilon ,\bar{w}\right] $ for all $\varepsilon \in \left( 0,\bar{w}\right) $.
    \item $\Pi \left( w,H_{i}\right) $ is linear over $\left[ \bar{w}-\varepsilon ,\bar{w}\right] $ for some maximal $\varepsilon \in \left( 0,\bar{w}\right)$.
\end{enumerate}
The rest of the proof rules out these possibilities. Define $b\left( w\right) \equiv \left(\ubar{U} -w\left( 1-c\right)\right) /\left( \mu -w\right) $, the inverse of $a\left(\cdot \right) $; and let $G_{i}\in \Delta \left( \left[\ubar{U},\bar{U}\right]
\right) $ be the equilibrium distribution of reservation values with density $g_{i}$. The next claim eliminates the first possibility.
\begin{claim}\label{strictconcavity}
$\Pi \left( \cdot ,H_{i}\right) $ is not strictly concave on some interval $\left[ \bar{w}-\varepsilon ,\bar{w}\right] $. Moreover, $\Pi \left( \cdot ,H_{i}\right) $ cannot
consist of an infinite sequence of linear segments on any neighborhood of $\bar{w}$.
\end{claim}
\begin{proof}
If  $\Pi \left( \cdot ;H_{i}\right) $ is strictly concave on some interval $\left[ \bar{w}-\varepsilon ,\bar{w}\right] $, it must be
increasing on $\left[ \bar{w}-\varepsilon ,\bar{w}\right] $. Thus, effective values $\left[ \bar{w}-\varepsilon ,\bar{w}\right] $ must be on the support of the equilibrium distribution $H_{i}$. Moreover, effective value $w\in \left[\bar{w}-\varepsilon ,\bar{w}\right] $ can be on the support of \textbf{only} the (pure-strategy) effective-value distribution with reservation value $%
b\left( w\right) $. Consequently, reservation
values in $\left[ b\left( \bar{w}\right) ,b\left( \bar{w}%
-\varepsilon \right) \right] $ must be on the support of $G_{i}$. Moreover,
the equilibrium condition that all reservation values in $\left[
b\left( \bar{w}\right) ,b\left( \bar{w}-\varepsilon \right) \right] $ yield
an expected payoff of $1/n$ implies that $\Pi \left( \cdot ;H_{i}\right) $
is necessarily strictly convex over $\left[ b\left( \bar{w}\right) ,b\left( 
\bar{w}-\varepsilon \right) \right] $. Therefore, effective values in $\left[ b\left( \bar{w}\right) ,b\left( \bar{w}-\varepsilon \right) %
\right] $ are realized \textit{only as reservation values} (rather than
posteriors). Furthermore, as $\Pi \left( \cdot ;H_{i}\right) $
is flat over $\left[ \bar{w},\ubar{U}\right] $ (if the interval is
nontrivial) and $\Pi \left( \ubar{U};H_{i}\right) =1/n$, the
payoff function must be flat at $1/n$ over $\left[ 
\ubar{U},b\left( \bar{w}\right) \right] $; otherwise, a payoff
exceeding $1/n$ can be attained by some reservation value in this interval.
The observations above suggest the following must hold.
\begin{enumerate}[label={(\roman*)},noitemsep,topsep=0pt]
    \item \label{concavpt1} $\left[ b\left( \bar{w}\right) ,b\left( \bar{w}-\varepsilon \right)\right] \in \supp\left(G_{i}\right) $, and $G_{i}\left( b\left( \bar{w}\right) \right) =0$.
    \item \label{concavpt2} Any reservation value $U\in \left[b\left(\bar{w}\right) ,b\left( \bar{w}-\varepsilon \right) \right] $ is realized by the (uniquely optimal) binary effective-value distribution supported on $\left\{a\left(U\right), U\right\} $, with respective weights $1-c/\left( 1-U\right) $ and $c/\left(1-U\right) $.
    \item \label{concavpt3} Reservation values outside of $\left[ b\left(\bar{w}\right),b\left(\bar{w}-\varepsilon \right) \right] $ assign zero measure to effective values in $\left[ \bar{w}-\varepsilon,\bar{w}\right] $.
\end{enumerate}
Using \ref{concavpt2} and $H_{i}\left( b\left( \bar{w}\right) \right)
=\Pi \left( b\left( \bar{w}\right) ;H_{i}\right)^{\frac{1}{n-1}}=\Pi \left( 
\ubar{U};H_{i}\right)^{\frac{1}{n-1}}=\left( n^{-1}\right)^{\frac{1}{%
n-1}}$, for reservation value $U\in \left[ b\left( \bar{w}\right) ,b\left( 
\bar{w}-\varepsilon \right) \right] $,
\[\label{eqm payoff 1} \tag{$A2$}H_{i}\left( U\right) =\left( \frac{1}{n}\right)^{\frac{1}{n-1}%
}+\int_{b\left( \bar{w}\right) }^{U}\frac{c}{1-u}dG_{i}\left( u\right) \text{ .} \]
Using \ref{concavpt3}, for reservation value $U\in \left[ b\left( \bar{w}\right)
,b\left( \bar{w}-\varepsilon \right) \right] $, we have 
\[\label{eqm payoff 2}\tag{$A3$}H_{i}\left( U\right) -H_{i}\left( a\left(U\right) \right) =G_{i}\left(
U\right) \text{ .} \]
We can combine three things--i. reservation values $U\in \left[ b\left( \bar{w}\right)
,b\left( \bar{w}-\varepsilon \right) \right] $ must yield payoff $1/n$, ii. $\Pi \left( w;H_{i}\right) =H_{i}\left( w\right)^{n-1}$ for all $w\in \left( 0,\bar{U}\right) $, and iii. (\ref{eqm payoff 2})--to obtain
\[\label{eqm payoff 3} \tag{$A4$}
\left( \frac{c}{1-U}\right) H_{i}\left( U\right)^{n-1}+\left( 1-\frac{c}{1-U%
}\right) \left( H_{i}\left( U\right) -G_{i}\left( U\right) \right)^{n-1}=%
\frac{1}{n}\text{ .}\]
Using (\ref{eqm payoff 1}), $H_{i}^{\prime }\left(U\right) =cg_{i}\left(
U\right) /\left( 1-U\right) $. Differentiating (\ref{eqm payoff 3}) with
respect to $U$ and substituting in for $H_i^{'}$, we obtain
\[\label{const payoff} \tag{$A5$} 0=c\left[ H_{i}\left( U\right)^{n-1}-\left( H_{i}\left( U\right)
-G_{i}\left( U\right) \right)^{n-1}\right] -\left( n-1\right) \left[ \left(
1-c-U\right)^{2}\left( H_{i}\left( U\right) -G_{i}\left( U\right) \right)
^{n-2}-c^{2}H_{i}\left( U\right)^{n-2}\right] g_{i}\left( U\right) \text{.}\]
Substituting $U=b\left(\bar{w}\right) $ into (\ref{const payoff}) gives $g_{i}\left( b\left( \bar{w}\right) \right) =0$. For any $U\in (b\left( \bar{w}\right) ,b\left( \bar{w}-\varepsilon \right) ]$, $G_{i}\left( U\right) >0$, so (\ref{const payoff}) implies $\left( 1-c-U\right)^{2}\left( H_{i}\left(U\right) -G_{i}\left( U\right) \right)^{n-2}-c^{2}H_{i}\left( U\right)
^{n-2}>0$. Rearranging (\ref{const payoff}) yields
\[\label{ratio} \tag{$A6$}
\frac{g_{i}\left( U\right) }{G_{i}\left( U\right) }=\frac{c/\left(
n-1\right) }{\left( 1-c-U\right)^{2}\left( H_{i}\left( U\right)
-G_{i}\left( U\right) \right)^{n-2}-c^{2}H_{i}\left( U\right)^{n-2}} \left[ \frac{H_{i}\left( U\right)^{n-1}-\left( H_{i}\left( U\right)
-G_{i}\left( U\right) \right)^{n-1}}{G_{i}\left( U\right) }\right] \equiv
\Upsilon \left( U\right) \text{ .}\]
By L'H\^{o}pital's rule, $\Upsilon \left( b\left( \bar{w}\right) \right)
=c/\left( \left( 1-\mu \right)^{2}-c^{2}\right) $. As $\Upsilon \left(
U\right) $ is continuous, it can be bounded over $[b\left( \bar{w}\right)
,b\left( \bar{w}-\varepsilon \right) ]$; define $\bar{\Upsilon}\equiv
\max_{\left\{ b\left( \bar{w}\right) ,b\left( \bar{w}-\varepsilon \right)
\right\}}\Upsilon \left( U\right) $. Integrating both sides of (\ref{ratio})
over an interval $\left[ U_{0},b\left( \bar{w}%
-\varepsilon \right) \right] $, we obtain
\[
\ln G_{i}\left( U_{0}\right) =\ln G_{i}\left( b\left( \bar{w}-\varepsilon
\right) \right) -\int_{U_{0}}^{b\left( \bar{w}-\varepsilon \right) }\Upsilon
\left( u\right) du\geq \ln G_{i}\left( b\left( \bar{w}-\varepsilon \right)
\right) -\bar{\Upsilon}\times \left( b\left( \bar{w}-\varepsilon \right) -%
\ubar{U}\right) \text{ .}\]
As the right-hand side of this is finite and $G_{i}$ has no atom at $b\left( \bar{w}\right) $, taking $U_{0}\rightarrow b\left( \bar{w}\right) $ yields a contradiction. 

We can adapt this proof to rule out equilibria with an infinite sequence
of decreasing ``kink'' points $\left\{ U_{k}\right\} $ such that i. $U_{k}\rightarrow b\left( \bar{w}\right) $, ii. $\Pi \left( \cdot
;H_{i}\right) $ is linear on each $\left[ U_{k+1},U_{k}\right] $,
and iii. $\Pi \left( \cdot ;H_{i}\right) $ is non-differentiable at each of
these kink points $U_{k}$.

First, (\ref{eqm payoff 1}) holds for all $U\in \left[ b\left( \bar{w}%
\right) ,U_{0}\right] $. Second, following a similar argument as above, (\ref%
{eqm payoff 2}) and hence (\ref{eqm payoff 3}) hold at all kink points. An
analogue of (\ref{const payoff}) can be obtained by subtracting (\ref{eqm
payoff 3}) for $U_{k}$ from that for $U_{k+1}$ and invoking the mean value
theorem as follows:
\[\begin{split}
    0 &\leq -\left( n-1\right) \left( \left( 1-c-U_{k}\right) \left(1-c-U_{k+1}\right) \left( H_{i}\left( U_{k}\right) -G_{i}\left( U_{k}\right)\right)^{n-2}-c^{2}H_{i}\left( U_{k+1}\right)^{n-2}+O\left( g_{i}\left(U_{k+1}\right) \left( U_{k}-U_{k+1}\right) \right) \right) g_{i}^{+}\left(U_{k+1}\right) \\
    &+c\left[ H_{i}\left( U_{k+1}\right)^{n-1}-\left( H_{i}\left(U_{k+1}\right) -G_{i}\left( U_{k+1}\right) \right)^{n-1}\right] \text{ ,}
\end{split}\]
where $g_{i}^{+}$ stands for the right-hand derivative of $G_{i}$. Dividing
both sides by $G_{i}\left( U_{k+1}\right) $, taking the limit (as $k\to \infty$) and rearranging, $\lim g_{i}^{+}\left( U_{k+1}\right) /G_{i}\left( U_{k+1}\right) $
is bounded. On the other hand, the linearity of $\Pi \left( \cdot
;H_{i}\right) $ over $\left[ U_{k+1},U_{k}\right] $ implies that $g\left(
U\right) $ is decreasing in $U$ over the interval.\footnote{%
The reason is as follows. $H_{i}\left( U\right) $ can be expressed as $%
\left( L_{0}+\sigma U\right)^{\frac{1}{n-1}}$ for some constants $L_{0}$
and $\sigma >0$. Therefore, $H_{i}^{\prime }\left( U\right) =\sigma \left(
L_{0}+\sigma U\right)^{-\frac{n-2}{n-1}}/\left( n-1\right) $. Equating this
with $H_{i}^{\prime }\left( U\right) =cg_{i}\left( U\right) /\left(
1-U\right) $ (recall (\ref{eqm payoff 1})), we have $g_{i}\left( U\right)
=\sigma \left( L_{0}+\sigma U\right)^{-\frac{n-2}{n-1}}\left( 1-U\right) /%
\left[ c\left( n-1\right) \right] $, which is clearly decreasing in $\sigma $%
.} Therefore, $g_{i}\left( U\right) /G_{i}\left( U\right) \leq
g_{i}^{+}\left( U_{k+1}\right) /G_{i}\left( U_{k+1}\right) $ for all $U\in
\left( U_{k+1},U_{k}\right) $. The boundedness of $\lim g_{i}^{+}\left(
U_{k+1}\right) /G_{i}\left( U_{k+1}\right) $ thus implies that $%
\lim_{U\rightarrow b\left( \bar{w}\right) }g_{i}\left( U\right) /G_{i}\left(
U\right) $ is bounded for all $\left[ b\left( \bar{w}\right) ,U_{0}\right] $. This is a contradiction, as argued above.
\end{proof}
We finish by ruling out the possibility that $\Pi \left( \cdot ;H_{i}\right) $ is
linear on some maximal $\left[ \bar{w}-\varepsilon ,\bar{w}\right] \subset \left[ 0,\ubar{U}\right]$.
\begin{claim}\label{nomaximalint}
For all $\varepsilon \in \left(0, \bar{w}\right)$ and $\delta \in \left(0, \bar{w} - \varepsilon\right)$, $\Pi \left( \cdot;H_{i}\right) $ is linear on $\left[ \bar{w}-\varepsilon -\delta
,\bar{w}\right] $.
\end{claim}
\begin{proof}
There are two cases to consider: i. $\bar{w}=\ubar{U}$ and ii. $\bar{w}<\ubar{U}$.

\textbf{Case i. $\bar{w}=\ubar{U}$:} Suppose $\Pi \left( \cdot ;H_{i}\right) $ is linear on $\left[\ubar{U}-\varepsilon ,\ubar{U}\right] $, where $\varepsilon \in \left( 0,\bar{w}\right) $ is maximal. We first explain why $\Pi \left( \cdot ;H_{i}\right) $ must be linear on $\left[\ubar{U}-\varepsilon ,b\left( \ubar{U}-\varepsilon \right) \right]$. Denote by $\sigma >0$ the slope of the linear segment of the payoff
function on $\left[ \ubar{U}-\varepsilon ,\ubar{U}\right] $. Suppose there is some $U^{\prime }\in (\ubar{U},b\left( \ubar{U}-\varepsilon \right) ]$ such that $\Pi \left( U^{\prime };H_{i}\right) >\Pi\left( \ubar{U};H_{i}\right) +\sigma \left( U^{\prime }-\ubar{U}\right) $. Then an inducible effective-value distribution with binary support $\left\{ a\left(U^{\prime }\right) ,U^{\prime }\right\} $ can achieve a payoff exceeding the equilibrium level $\Pi \left( \ubar{U}%
;H_{i}\right) $, a contradiction. 

Next, suppose there is some $U^{\prime }\in (\ubar{U},b\left( \ubar{U}-\varepsilon \right) ]$ such that $\Pi\left( U^{\prime };H_{i}\right) <\Pi \left( \ubar{U};H_{i}\right)
+\sigma \left( U^{\prime }-\ubar{U}\right) $. Effective value $U^{\prime }$ cannot deliver the equilibrium payoff and hence cannot lie on the support of the equilibrium effective-value distribution $H_{i}$. As there are no interior atoms, this implies $U^{\prime }\geq \hat{U}$ and $\Pi \left(w;H_{i}\right) <\Pi \left( \ubar{U};H_{i}\right) +\sigma \left( w-%
\ubar{U}\right)$ for all $w\in \left[ U^{\prime },\bar{U}\right)$.
Moreover, because of the maximality of $\varepsilon $ and the concavity of $%
\Pi _{\ubar{U}}\left( \cdot ;H_{i}\right) $, we have $\Pi \left(
w;H_{i}\right) <\Pi \left( \ubar{U};H_{i}\right) +\sigma \left( w-%
\ubar{U}\right) $ for $w\in \left( 0,\ubar{U}-\varepsilon \right) $%
. Therefore, for reservation values less than $b\left( \ubar{U}-\varepsilon \right) $, it is strictly suboptimal to assign positive weight
to effective values in $\left( 0,\ubar{U}-\varepsilon
\right) $. This, plus $b\left(\ubar{U}-\varepsilon \right) > \hat{U}$, implies effective values in the interval $\left( 0,\ubar{U}-\varepsilon \right) 
$ cannot lie on the support of the equilibrium distribution $H_{i}$. The
flatness of $\Pi \left( \cdot ;H_{i}\right) $ over interval $\left( 0,%
\ubar{U}-\varepsilon \right) $ contradicts Lemma \ref{concavity}.

Next, we explain that in equilibrium, the effective-value distribution on $\left[  \ubar{U}-\varepsilon ,b\left(  \ubar{U}-\varepsilon \right) \right] $ is "self-contained:" any reservation value $U\in \supp\left\{ G_{i}\right\} \cap \left[  \ubar{U},b\left( \ubar{U}-\varepsilon \right) \right] $ is realized by an effective-value
distribution supported on $\left[  \ubar{U}-\varepsilon ,U\right] $ and
any reservation value $U\in \supp\left\{ G_{i}\right\} $ exceeding $b\left( 
 \ubar{U}-\varepsilon \right) $ must be realized by an effective-value
distribution that puts no mass in $\left[\ubar{U}-\varepsilon ,b\left(  \ubar{U}-\varepsilon \right) \right] $. The first observation follows from the maximality of $\varepsilon $; specifically, for
all $w\in \left( 0, \ubar{U}-\varepsilon \right) $, $\Pi \left(
w;H_{i}\right) <\Pi \left(  \ubar{U};H_{i}\right) +\sigma \left( w-%
 \ubar{U}\right) $. To see the second observation, let $U\in \supp\left\{
G_{i}\right\} $ be some reservation value exceeding $b\left(  \ubar{U}%
-\varepsilon \right) $, so that $a\left(U\right) < \ubar{U}-\varepsilon 
$. If its corresponding effective-value distribution assigns positive mass
to $\left[  \ubar{U}-\varepsilon ,b\left(  \ubar{U}%
-\varepsilon \right) \right] $, then it must also assign positive mass to $[0,a\left(U\right) )$. $\varepsilon $'s maximality implies $\Pi\left( a\left(U\right) ;H_{i}\right) >\lambda \Pi \left( w_{0};H\right) +\left( 1-\lambda \right) \Pi \left( w_{1};H_{i}\right) $ for any combination
of $w_{0}$ and $w_{1}$ such that $w_{1}\geq  \ubar{U}-\varepsilon $, $w_{0}<a\left(U\right) $ and $\lambda w_{0}+\left( 1-\lambda \right) w_{1}=a\left(U\right) $. Consequently, this cannot be an optimal effective value distribution given $U$.

These facts imply that $H_{i}\left( w\right) =\left( n^{-1}+\sigma \left( w- \ubar{U}\right)
\right)^{\frac{1}{n-1}}$ for $w\in \left[ \ubar{U}-\varepsilon
,b\left(  \ubar{U}-\varepsilon \right) \right] $ and that
\[H_{i}^{\#}\left( w\right) \equiv \frac{H_{i}\left( w\right) -H_{i}\left( 
 \ubar{U}-\varepsilon \right) }{H_{i}\left( b\left(  \ubar{U}%
-\varepsilon \right) \right) -H_{i}\left(  \ubar{U}-\varepsilon \right) }
\text{ ,}\]
the conditional distribution for $w\in \left[  \ubar{U}-\varepsilon
,b\left(  \ubar{U}-\varepsilon \right) \right] $, is inducible. The rest of the proof derives a contradiction by showing that the mean condition $\int_{\ubar{U} - \varepsilon}^{b\left(  \ubar{U}-\varepsilon\right) }wdH_{i}^{\#}\left( w\right) = \ubar{U}$ is impossible to satisfy.

Our first step in showing the unattainability of the mean condition is to derive an upper bound on the slope parameter $\sigma $. On one hand, the concavity of $\Pi \left( \cdot ;H_{i}\right) $ over $(0, \ubar{U}]$ implies that $\sigma $ cannot exceed $\left(1/n-\lim_{w\rightarrow 0^{+}}\Pi \left( w;H_{i}\right) \right)/\ubar{U}$. On the other hand, as $\Pi \left( b\left(  \ubar{U}%
-\varepsilon \right) ;H_{i}\right) $ is at most $\lim_{w\rightarrow \bar{U}%
^{-}}\Pi \left( w;H_{i}\right) $, $\sigma $ cannot exceed $\left(
\lim_{w\rightarrow \bar{U}}\Pi \left( w;H_{i}\right) -n^{-1}\right) /\left(
b\left(  \ubar{U}-\varepsilon \right) - \ubar{U}\right) $ either.
Recall $\lim_{w\rightarrow 0^{+}}\Pi \left( w;H_{i}\right) =\left( \alpha
\left( 1-\mu \right) \right)^{n-1}$ and $\lim_{w\rightarrow \bar{U}^{-}}\Pi
\left( w;H_{i}\right) =\left( 1-\alpha \mu \right)^{n-1}$. An upper bound $%
\bar{\sigma}\left( \varepsilon \right) $ on the slope $\sigma $ is thus%
\[\label{sigmabar}\tag{$A7$}\bar{\sigma}\left( \varepsilon \right) \equiv \min \left\{\frac{\frac{1}{n}-\left( \alpha \left( 1-\mu \right) \right)^{n-1}}{\mu -c},\frac{%
\left( 1-\alpha \mu \right)^{n-1}-\frac{1}{n}}{1-\mu}\frac{c+\varepsilon }{%
\varepsilon }\right\} = \begin{cases}
\frac{\frac{1}{n}-\left( \alpha \left( 1-\mu \right) \right)^{n-1}}{\mu -c} \ &\text{if} \ \varepsilon \leq  \ubar{U}-a\left(\hat{U}\right) \\ 
\frac{\left( 1-\alpha \mu \right)^{n-1}-\frac{1}{n}}{1-\mu}\frac{c+\varepsilon 
}{\varepsilon} \ & \text{if} \ \varepsilon > \ubar{U}-a\left(\hat{U}\right)
\end{cases}  \text{ .}\]
Observe that $\bar{\sigma}\left( \varepsilon \right) $ is nonincreasing in $\varepsilon $.\footnote{%
To see that the cutoff value of $\varepsilon $ is related to $\hat{U}$, note
that $\hat{U}$ is the unique reservation value such that points $\left(
0,\lim_{w\rightarrow 0}\Pi \left( w;H_{i}\right) \right) $, $\left( 
 \ubar{U},\Pi \left(  \ubar{U},H_{i}\right) \right) $, and $\left( 
\hat{U},\lim_{w\rightarrow \bar{U}}\Pi \left( w;H_{i}\right) \right) $ all
lie on the same straight line on the graph of $\Pi \left( \cdot
;H_{i}\right) $.} The mean condition for $H_{i}^{\#}$'s inducibility can be expressed as $\Psi \left( \sigma ,\varepsilon \right) =0$, where
\[
\Psi \left( \sigma ,\varepsilon \right) \equiv \int_{\ubar{U}%
-\varepsilon }^{b\left( \ubar{U}-\varepsilon \right) }H_{i}^{\#}\left(
w\right) dw-\left( b\left( \ubar{U}-\varepsilon \right) -\ubar{U}%
\right) \text{ .}\]
Directly, $\Psi\left( 0,\varepsilon \right) =0$ and 
\[
\frac{d\Psi \left( \sigma ,\varepsilon \right) }{d\sigma }=\frac{\sigma\varepsilon^{2}}{n-1}\left[\left(\frac{1}{n}-\sigma \varepsilon\right)^{-\frac{n-2}{n-1}}-\left(\frac{1-\mu}{c+\varepsilon}\right)^{2}\left(\frac{1}{n}+\sigma \left(\frac{1-\mu }{c+\varepsilon}\varepsilon\right)\right)^{-\frac{n-2}{n-1}}\right] \text{ .}
\]
As the bracketed term is increasing in $\sigma $, it can change sign at most
once as $\sigma $ increases, and the change can only be from negative to
positive. Thus, either $\Psi \left(\sigma ,\varepsilon \right) $ is
increasing in $\sigma $, or is a U-shaped function in $\sigma $ (fixing $\varepsilon $). In the former case, $\Psi \left( \sigma ,\varepsilon \right) 
$ clearly has no non-zero solution. Therefore, we focus on the latter case.
It suffices to show that $\Psi \left( \bar{\sigma}\left( \varepsilon \right)
,\varepsilon \right) <0$ for all $\varepsilon \in (0,\ubar{U})$.

Consider the case $\varepsilon \leq \ubar{U}-a\left(\hat{U}\right) $.
Note that $\bar{\sigma}\left( \varepsilon \right) $ is independent of $\varepsilon $ (see (\ref{sigmabar})). $\Psi \left( \bar{\sigma}\left( 0\right) ,0\right) =0$, and 
\[
\frac{d\Psi \left( \bar{\sigma}\left( \varepsilon \right) ,\varepsilon
\right) }{d\varepsilon }=\frac{\sigma^{2}\varepsilon }{n-1}\left[\left(\frac{1}{n}-\sigma \varepsilon \right)^{-\frac{n-2}{n-1}}-\frac{\left(1-\mu \right)^{2}c}{\left( c+\varepsilon \right)^{3}}\left( \frac{1}{n}+\sigma \left( 1-\mu \right) \frac{\varepsilon }{c+\varepsilon }\right)^{-\frac{n-2}{n-1}}\right] \text{ .}\]
The bracketed term is increasing in $\varepsilon $, so it can
change sign at most once as $\varepsilon $ increases (from negative to positive). Consequently, it suffices to check $\Psi \left( \bar{\sigma}\left( \varepsilon^{*} \right) ,\varepsilon^{*} \right) <0$ where $\varepsilon^{*} \equiv \ubar{U}-a\left(\hat{U}\right) $. Directly, 
\[\Psi \left( \bar{\sigma}\left(\varepsilon^{*} \right)
,\varepsilon^{*} \right) =\int_{a\left(\hat{U}\right) }^{\hat{U}}H_{i}^{\#}\left( w\right) dw-\left( \hat{U}-\ubar{U}\right) < \int_{a\left(\hat{U}\right) }^{\hat{U}}\frac{H_{i}\left( w\right) -\alpha
\left( 1-\mu \right) }{\left( 1-\alpha \mu \right) -\alpha \left( 1-\mu
\right) } dw-\left( \hat{U}-\ubar{U}\right) = 0 \text{ .}\]

Finally, consider the case $\varepsilon \in \left( \ubar{U}-a\left(\hat{U}%
\right) ,\ubar{U}\right) $. Using (\ref{sigmabar}), 
\[\frac{d\Psi \left( \bar{\sigma}\left( \varepsilon \right) ,\varepsilon
\right) }{d\varepsilon }=\left( n^{-1}-\frac{\left( 1-\alpha \mu \right)
^{n-1}-n^{-1}}{1-\mu }\left( c+\varepsilon \right) \right)^{-\frac{n-2}{n-1}%
}\frac{n\left( \left( 1-\alpha \mu \right)^{n-1}-n^{-1}\right)^{2}\left(
c+\varepsilon \right) }{\left( n-1\right) \left( 1-\mu \right)^{2}}>0 \text{ .}\]
 At $\varepsilon = \varepsilon^{\dagger} \equiv \ubar{U}$ and $\sigma =\bar{\sigma}\left( \ubar{U}\right) $,
\[\Psi \left( \bar{\sigma}\left( \varepsilon^{\dagger} \right) ,\varepsilon^{\dagger} \right)
=\int_{0}^{\bar{U}}H_{i}^{\#}\left( w\right) dw-\left( \bar{U}-\ubar{U}%
\right) <\int_{0}^{\bar{U}}\frac{H_{i}\left( w\right) -\alpha
\left( 1-\mu \right) }{\left( 1-\alpha \mu \right) -\alpha \left( 1-\mu
\right) } dw-\left( \bar{U}-%
\ubar{U}\right) =0 \text{ ,}\]
since $n^{-1}-\bar{%
\sigma}\left( \ubar{U}\right) \ubar{U}>\left( \alpha \left( 1-\mu
\right) \right)^{n-1}$ (recall the
definition of $\bar{\sigma}\left( \cdot \right) $, which implies $\bar{\sigma}\left( \ubar{U}\right) <\left( n^{-1}-\left( \alpha
\left( 1-\mu \right) \right)^{n-1}\right) /\ubar{U}$).

\textbf{Case ii. $\bar{w} < \ubar{U}$:} We begin the analysis of this case with a number of observations. First, as $\Pi\left( \cdot ;H_{i}\right) $ is flat over $\left[ \bar{w},\ubar{U}\right] $, it must also be flat over $\left[\ubar{U},b\left( \bar{w}%
\right) \right] $; otherwise, a firm could achieve a payoff exceeding $\Pi \left(\ubar{U};H_{i}\right) $ (by adopting a reservation value between $\ubar{U}$ and $b\left( \bar{w}\right) $). Second, the modified concavification approach tells us reservation values in the interval $\left[ b\left( \bar{w}\right) ,b\left( \bar{w}-\varepsilon \right) \right] $ must be on the support of the equilibrium
reservation-value distribution $G_{i}$. Third, (in equilibrium) any reservation value $U$ in the interval $\left[b\left( \bar{w}\right) ,b\left( \bar{w}-\varepsilon \right) \right] $ can
only be realized using an effective-value distribution with the following
form: a weight $c/\left( 1-U\right) $ is put on $U$, and the complementary
weight $1-c/\left( 1-U\right) $ to values in the interval $\left[ \bar{w}%
-\varepsilon ,\bar{w}\right] $ in such a way that the conditional mean is $%
a\left(U\right) $. \footnote{To see this, it suffices to note that no weight should be put on any values
in the interval $\left( b\left( \bar{w}\right) ,U\right) $. Suppose not. Then there must be some $w\in \left( b\left( \bar{w}\right),U\right) $ such that on the graph of $\Pi \left( \cdot ;H\right) $, $\left(w,\Pi \left( w;H\right) \right) $ lies on the straight line extending from that connecting $\left( \bar{w}-\varepsilon,\Pi \left( w_{1};H\right) \right) $ and $%
\left( \bar{w},\Pi \left( \bar{w};H\right) \right) $. It is then apparent that reservation value $w$ yields a payoff strictly exceeding $\Pi \left(\ubar{U};H\right) $, a contradiction.} Finally, as the optimal effective-value distribution conditional on reservation value $U\in \left[b\left( \bar{w}\right) ,b\left( \bar{w}-\varepsilon \right) \right] $ must
yield the equilibrium payoff $1/n$, the payoff function $\Pi \left( \cdot
;H_{i}\right) $ on $w\in \left[ \bar{w}-\varepsilon ,b\left( 
\bar{w}-\varepsilon \right) \right] $ necessarily takes the following form:
\[\Pi \left( w;H_{i}\right) = \begin{cases}
\frac{1}{n}+\sigma \left( w-\bar{w}\right), \quad  &\text{if} \quad w \in \left[\bar{w}%
-\varepsilon ,\bar{w}\right) \\ 
\frac{1}{n}, \quad &\text{if} \quad w\in \left[ \bar{w},b\left( \bar{w}\right) \right] \\ 
\frac{1}{n}-\sigma \frac{\left( \mu -c\right) -\left( 1-c\right) \bar{w}-\left(
\mu -\bar{w}\right) w}{c}, \quad  &\text{if} \quad w\in (b\left( \bar{w}\right)
,b\left( \bar{w}-\varepsilon \right) ]
\end{cases} \text{ ,}  \label{payoff linear flat segment} \tag{$A8$}\]
for some slope parameter $\sigma >0$.

The third observation above tells us that on $\left[ b\left( \bar{w}%
\right) ,b\left( \bar{w}-\varepsilon \right) \right] $, $H_{i}$ and $G_{i}$
are related by (\ref{eqm payoff 1}). Consequently, $H_{i}^{\prime }\left(
w\right) =\left[ c/\left( 1-w\right) \right] \times g_{i}\left( w\right)$. 
Moreover, using (\ref{payoff linear flat segment}), $d\Pi \left( w;H_{i}\right) /dw=\sigma \left( \mu -%
\bar{w}\right) /c$. As $\Pi \left( w;H_{i}\right) =H_{i}\left( w\right)
^{n-1}$, we have $g_{i}\left( w\right) =\left[ \sigma \left( \mu -\bar{w}%
\right)/\left(c^{2}\left( n-1\right)\right) \right] \times \left[ \left( 1-w\right)
/H_{i}\left( w\right) ^{n-2}\right] $, so
\[g_{i}\left( b\left( \bar{w}\right) \right) =\frac{\sigma n^{\frac{n-2}{n-1}}}{c\left( n-1\right) }\left( 1-\bar{w}\right) \text{ .}
\label{boundary density} \tag{$A9$}
\]
Below, we derive a counterpart of (\ref{eqm payoff 2}) for this case. Note that because of the locally linearity of the payoff function, we can no
longer guarantee that reservation value $U\in \left[ b\left( \bar{w}\right)
,b\left( \bar{w}-\varepsilon \right) \right] $ is generated uniquely by the
binary effective-value distribution supported on $\left\{ a\left(U\right)
,U\right\} $. Nonetheless, we can still rely on the third observation above
to derive an bound on $G_{i}\left( U\right) $. Let $l$ be some positive real
number. For any reservation value $u\in \left[ b\left( \bar{w}\right) ,U\right] $, its corresponding effective-value distribution assigns a mass of $%
1-c/\left( 1-u\right) $ to values in the interval $\left[ \bar{w}%
-\varepsilon ,\bar{w}\right] $ in such a way that the conditional mean is $%
a\left(u\right) $. Therefore, by Bayes' rule, the mass that it assigns to
values outside of the interval $\left[ a\left(u\right) -l\left( \bar{w}%
-a\left(u\right) \right) ,\bar{w}\right] $ cannot exceed $\left( 1-c/\left(
1-u\right) \right) \times \left( \bar{w}-a\left(u\right) \right) /\left[
\left( 1+l\right) \left( \bar{w}-a\left(U\right) \right) \right] $. As a
result, for any $U\in \left[ b\left( \bar{w}\right) ,b\left( \bar{w}%
-\varepsilon \right) \right] $,%
\[
G_{i}\left( U\right) \leq H_{i}\left( U\right) -H_{i}\left( a\left(U\right)
-l\left( \bar{w}-a\left(U\right) \right) \right) +\int_{b\left( \bar{w}%
\right) }^{U}\left( 1-\frac{c}{1-u}\right) \frac{\bar{w}-a\left(u\right) }{%
\left( 1+l\right) \left( \bar{w}-a\left(U\right) \right) }dG_{i}\left(
u\right) \text{ .}\]
Using $H_{i}\left( \cdot \right) =\Pi \left( \cdot
;H_{i}\right) ^{\frac{1}{n-1}}$ and the payoff specification of (\ref{payoff
linear flat segment}), the inequality above can be rewritten as
\[\begin{split}
&G_{i}\left( U\right) -\int_{b\left( \bar{w}\right) }^{U}\left( 1-\frac{c}{1-u}\right) \frac{\bar{w}-a\left(u\right) }{\left( 1+l\right) \left( \bar{w}-a\left(U\right) \right) }dG_{i}\left( u\right)\\
&\leq \left( n^{-1}-\frac{\sigma }{c}\left[ \left( \mu -c\right) -\left(1-c\right) \bar{w}-\left( \mu -\bar{w}\right) U\right] \right) ^{\frac{1}{n-1}}-\left( n^{-1}-\sigma \left( 1+l\right) \left( \bar{w}-a\left(U\right)\right) \right) ^{\frac{1}{n-1}}
\end{split}\text{ .}  \label{limit ineq}\tag{$A10$}\]
Dividing the left-hand side of (\ref{limit ineq}) by $U-b\left( \bar{w}%
\right) $ and taking the limit $U\rightarrow b\left( \bar{w}\right) $ gives:%
\[\left. \frac{d}{dU}\left( G_{i}\left( U\right) -\int_{b\left( \bar{w}%
\right) }^{U}\left( 1-\frac{c}{1-u}\right) \frac{\bar{w}-a\left(u\right) }{%
\left( 1+l\right) \left( \bar{w}-a\left(U\right) \right) }dG_{i}\left(
u\right) \right) \right|_{U=b\left( \bar{w}\right) } = g_{i}\left( b\left( \bar{w}\right) \right) \left( 1-\frac{1}{2\left(1+l\right) }\frac{1-\mu }{1-\bar{w}}\right) \text{ ,}\]
via L'H\^{o}pital's rule and direct substitution of $b\left( \cdot \right) $. Doing likewise to the right-hand side of (\ref{limit ineq}) yields
\[\frac{\sigma }{c}\frac{n^{\frac{n-2}{n-1}}}{n-1}\left( \mu -\bar{w}\right) \left( 1+\left( 1+l\right) \frac{\mu -\bar{w}}{1-\mu }\right) \text{ ,}\]
again via L'H\^{o}pital's rule. These computations, (\ref{boundary density}), and (\ref{limit ineq}) imply $\bar{w}\leq \mu -\left( 1-\mu \right)\frac{\sqrt{l+1/2}}{1+l}$. This reduces to $\bar{w}\leq \left(1+1/\sqrt{2}\right) \mu -1/\sqrt{2}$ as $l\rightarrow 0$. Because $\bar{w}\geq 0$ this is impossible whenever $\mu \leq \sqrt{2}-1$.

As $\bar{\mu}>\sqrt{2}-1$ for $n\geq 4$, it remains to consider the cases of 
$n=2$ and $3$. To this end, note that the local piecewise linearity of the
payoff function (\ref{payoff linear flat segment}) implies that the
effective-value distribution conditional on the interval $\left[ \bar{w}%
-\varepsilon ,b\left( \bar{w}-\varepsilon \right) \right] $, denoted by $%
H^{\#}$, is given by 
\[H^{\#}\left( w\right) \equiv \begin{cases}
\left[ \left(\frac{1}{n}+\sigma \left( w-\bar{w}\right) \right) ^{\frac{1}{n-1}%
}-\left(\frac{1}{n}-\sigma \varepsilon \right) ^{\frac{1}{n-1}}\right]/\bar{H}, \quad &\text{if} \quad w\in \lbrack \bar{w}-\varepsilon ,\bar{w}) \\
\left[ n^{-\frac{1}{n-1}}-\left(\frac{1}{n}+\sigma \varepsilon \right) ^{\frac{1%
}{n-1}}\right]/\bar{H}, \quad &\text{if} \quad w\in \left[\bar{w},b\left( 
\bar{w}\right) \right]\\ 
\left[\left(\frac{1}{n}-c^{-1}\sigma \left( \left( \mu -c\right) -\left(
1-c\right) \bar{w}-w\left( \mu -\bar{w}\right) \right) \right) ^{\frac{1}{n-1%
}}-\left(\frac{1}{n}+\sigma \varepsilon \right) ^{\frac{1}{n-1}}\right]/\bar{H}, \quad &\text{if} \quad w\in (b\left( \bar{w}\right) ,b\left( \bar{w}-\varepsilon \right) ]
\end{cases} \text{ ,}\]
where
\[
\bar{H}\equiv \left[\frac{1}{n}-c^{-1}\sigma \left( \left( \mu -c\right) -\left(
1-c\right) \bar{w}-b\left( \bar{w}-\varepsilon \right) \left( \mu -\bar{w}%
\right) \right) \right] ^{\frac{1}{n-1}}-\left(\frac{1}{n}+\sigma \varepsilon
\right) ^{\frac{1}{n-1}}\text{ .}\]
$H^{\#}$ must be inducible. Necessarily, $\int_{\bar{w}-\varepsilon }^{b\left( \bar{w}-\varepsilon\right) }H^{\#}\left( w\right) dw-\left[ b\left( \bar{w}-\varepsilon \right)
-\ubar{U}\right] =0$, or $\Upsilon _{n}\left( \sigma
,\varepsilon \right) =0$, where 
\[\begin{split}
    \Upsilon _{n}\left( \sigma ,\varepsilon ,c\right)  \equiv &\frac{n-1}{n}\left[
\left( \mu -c-\bar{w}\right) \frac{1}{n}^{\frac{n}{n-1}}-\left( \mu -\bar{w}%
\right) \left( \frac{1}{n}-\sigma \varepsilon \right) ^{\frac{n}{n-1}%
}+c\left( \frac{1}{n}+\sigma \varepsilon \frac{1-\mu }{\mu -\bar{w}%
+\varepsilon }\right) ^{\frac{n}{n-1}}\right]  \\
&+\sigma \frac{1}{n}^{\frac{1}{n-1}}\left( 1-\bar{w}\right) \left( \mu -c-%
\bar{w}\right) -\sigma \left( \mu -\bar{w}\right) \left( \left( 1-\mu
\right) \frac{\mu -w-c+\varepsilon }{\mu -w+\varepsilon }\right) \left( 
\frac{1}{n}+\sigma \varepsilon \frac{1-\mu }{\mu -\bar{w}+\varepsilon }%
\right) ^{\frac{1}{n-1}} \\
&-\sigma \left( \mu -\bar{w}\right) \left( \mu -c-\left( \bar{w}%
-\varepsilon \right) \right) \left( \frac{1}{n}-\sigma \varepsilon \right) ^{%
\frac{1}{n-1}}
\end{split} \text{ .}\]

$\Upsilon _{n}\left( \sigma ,\varepsilon ,c\right) $ is
linear in $c$, so it suffices to show that $\Upsilon _{n}\left( \sigma
,\varepsilon ,c\right) <0$ at both $c=0$ and $c=\mu -\bar{w}$ for all
possible $\sigma $ and $\varepsilon $. $c=\mu -\bar{w}$
reduces to the case $\bar{w}=\ubar{U}$ studied above. It remains to
analyze $c=0$. For $n=2$, 
\[\Upsilon _{2}\left( \sigma ,\varepsilon ,0\right) =-\frac{1}{2}\sigma
^{2}\varepsilon \left( \frac{\mu -\bar{w}}{\mu -w+\varepsilon }\right) \left[
2\left( 1-\mu \right) ^{2}-\left( \mu -\bar{w}+\varepsilon \right) \left(
2\left( \mu -\bar{w}\right) +\varepsilon \right) \right] \text{ .}\]
As $\varepsilon \leq \bar{w}$ and $\bar{w}>0$, the term in brackets is positive. For $n=3$,%
\[\Upsilon _{3}\left( \sigma ,\varepsilon ,0\right) =\frac{\mu -\bar{w}}{\sqrt{%
3}}\left[\frac{2}{9}-\left( \frac{2}{9}+\sigma \left( \mu -\bar{w}\right) +\frac{1}{3}\sigma \varepsilon \right)\sqrt{1-3\sigma \varepsilon}+\sigma \left(1-\bar{w}\right) -\sigma \left(1-\mu\right)
\left(1+3\sigma \varepsilon \frac{1-\mu }{\mu -\bar{w}+\varepsilon}\right)^{\frac{1}{2}}\right] \text{ .}\]
The term in brackets is decreasing in $\bar{w}$. Moreover, define $x\equiv \sigma \varepsilon $ to be the vertical increment of the
linear segment of the payoff function $\Pi \left( \cdot ;H_{i}\right)$. $x$ cannot exceed $1/3-\left(\tilde{\alpha}\left( 1-\mu \right) \right)^{2}$, where $\tilde{\alpha}$ is the nontrivial solution of (\ref{alpha U_hat no u0 middle mean 1}) when $n=3$. Denote this upper bound by $\bar{x}\left( \mu \right) $. It suffices to show that for all $x\leq \bar{x}\left( \mu \right) $ and $\sigma \geq x/\mu $,
\[\hat{\Upsilon}_{3}\left( \sigma ,x\right) \equiv \frac{2}{9}+\sigma-x-\left( \frac{2}{9}+\sigma\mu-\frac{2}{3}x\right)\sqrt{1-3x}-\sigma \left(1-\mu \right) \left( 1+3x\frac{1-\mu }{\mu}\right)^{\frac{1}{2}}<0\text{.}\]
The result follows from the fact that $d\hat{\Upsilon}_{3}\left( \sigma
,x\right) /d\sigma <0$ and $\hat{\Upsilon}_{3}\left( x/\mu,x\right) <0$. \end{proof}
Finally, (\ref{linear payoff no u0}) and (\ref{payoff in H}) allow us to recover the equilibrium distribution over effective values.
\begin{claim}\label{eqdistexplicit}
The equilibrium distribution over effective values is
\[H_{i}\left( w\right) = \begin{cases}
\alpha \left( 1-\mu \right), \quad &\text{if} \quad w = 0\\
\left(\left( \alpha \left( 1-\mu \right) \right)^{n-1}+\frac{\left(
1-\alpha \mu \right)^{n-1}-\left( \alpha \left( 1-\mu \right) \right)^{n-1}}{\hat{U}}w\right)^{\frac{1}{n-1}}, \quad &\text{if} \quad \text{if }w\in \left(0,\hat{U}\right]
\\ 
1-\alpha \mu, \quad &\text{if} \quad w\in (\hat{U},\bar{U}) \\ 
1, \quad &\text{if} \quad w=\bar{U}
\end{cases} \text{ .}  \label{linear H no u0} \tag{$A11$}
\]
\end{claim}
\end{proof}

\subsection{Proof of Lemma \protect\ref{binary dist no u0}}
\begin{proof}
Because of the continuity of $H_{i}\left( w\right) $ over $\left( 0,\hat{U}\right] $, it suffices to show that in this region, there is a mixed strategy that can match its density, denoted by $h_{i}$:
\[h_{i}\left( w\right) =\frac{\left( 1-\alpha \mu \right)^{n-1}-\left( \alpha
\left( 1-\mu \right) \right)^{n-1}}{\left( n-1\right) \hat{U}}\left( \left(
\alpha \left( 1-\mu \right) \right)^{n-1}+\frac{\left( 1-\alpha \mu \right)
^{n-1}-\left( \alpha \left( 1-\mu \right) \right)^{n-1}}{\hat{U}}w\right)
^{-\frac{n-2}{n-1}}\text{.}  \label{h density}
\tag{$A12$}\]
To this end, define a mapping $\gamma\colon \left[ \ubar{U},\hat{U}\right]
\rightarrow \left[ 0,\ubar{U}\right] $ by $K\left( U\right) =K\left( \gamma\left( U\right) \right) $, where $K\colon\left[0,\bar{U}\right]
\rightarrow \mathbb{R}$ is
\[\begin{split}
    K\left( w\right) &\equiv \left( \left( \alpha \left( 1-\mu \right) \right)^{n-1}+\left( \left( 1-\alpha \mu \right)^{n-1}-\left( \alpha \left( 1-\mu\right) \right)^{n-1}\right) \frac{w}{\hat{U}}\right)^{\frac{1}{n-1}} \\
    &\times \left( n\ubar{U}+\left( n-1\right) \left( \frac{\hat{U}}{\left( 1-\alpha \mu \right)^{n-1}-\left( \alpha \left( 1-\mu \right)\right)^{n-1}}\right) \left( \alpha \left( 1-\mu \right) \right)^{n-1}-w\right)\end{split} \text{ ,}\]
and parameters $\alpha $ and $\hat{U}$ are as specified in Lemma \ref{linear structure no u0}. For each $U\in \left[ \ubar{U},\hat{U}\right] $, let $F\left(\cdot ; U\right)$ be a binary distribution with support $\left\{\gamma\left( U\right),U\right\} $ and mean $\ubar{U}$; and let $F\left(\cdot ; \bar{U}\right)$ be the binary distribution with support $\left\{0,\bar{U}\right\} $ and mean $\ubar{U}$. Moreover, let $G$ be a reservation-value distribution that has an atom $\alpha \in \left[ 0,1\right] $ at $\bar{U}$ and a density for $U\in \left[ \ubar{U},\hat{U}\right] $ as follows: \[g\left( U\right) \equiv \frac{\left( 1-\alpha \mu \right)^{n-1}-\left(\alpha \left( 1-\mu \right) \right)^{n-1}}{\left( n-1\right) \hat{U}}\left(\left( \alpha \left( 1-\mu \right) \right)^{n-1}+\frac{\left( 1-\alpha \mu\right)^{n-1}-\left( \alpha \left( 1-\mu \right) \right)^{n-1}}{\hat{U}}U\right)^{-\frac{n-2}{n-1}}\frac{U-\gamma\left(U\right) }{\ubar{U}-\gamma\left(U\right) } \text{ .}\]
We claim that the mixed strategy \[\tag{$A13$}\label{mixedstrat42}\left( G,\left\{F\left(\cdot ; U\right) \colon U\in \left[ \ubar{U},\hat{U}\right] \cup \left\{ \bar{U}\right\}\right\}\right) \text{ ,} \] generates the effective-value distribution $H_i$ (defined in (\ref{linear H no u0})).

First, we show that effective-value distribution $F\left(\cdot ; U\right)$ is inducible for each $U\in \left[ \ubar{U},\hat{U}\right] \cup \left\{ \bar{U}\right\} $. Note that the mapping $\gamma$ is well-defined: a direct computation reveals that $K\left( w\right) $ is strictly concave with a peak at $\ubar{U}$ and that $K\left( 0\right) =K\left( \hat{U}\right) $.
We need to show that $\gamma\left( U\right) \leq a\left(U\right) \equiv \frac{\mu -c-\mu U}{1-c-U}$ (defined in \ref{baradef}). Furthermore, note that because $a\left(\ubar{U}\right) = \gamma\left( \ubar{U}\right) $, $\gamma\left( \hat{U}\right) =0=a\left(\bar{U}\right) \leq a\left(\hat{U}\right) $, and $a\left( U\right) $ is decreasing and strictly concave, it suffices to show that $\gamma\left( U\right) $ is convex. To this end, we adopt a change of variable: let $v=U-\ubar{U}$, and $d\left( v\right) = \ubar{U}-\gamma\left( \ubar{U}+v\right) $.\footnote{Intuitively, $d\left( \cdot \right) $ represents the "reflection" of $U$
about $\ubar{U}$ according to the function $K$.} The implicit
definition of $\gamma$ implies $K\left( \ubar{U}+v\right) =K\left(\ubar{U}-d\left( v\right) \right) $; or equivalently,
\[\left( M+v\right)^{\frac{1}{n-1}}\left( \left( n-1\right) M-v\right) =\left( M-d\left( v\right) \right)^{\frac{1}{n-1}}\left( \left( n-1\right)M+d\left( v\right) \right) \text{ ,}  \label{d def}\tag{$A14$}\]
where $M\equiv \ubar{U}+\frac{\hat{U}\left( \alpha \left( 1-\mu
\right) \right)^{n-1}}{\left( 1-\alpha \mu \right)^{n-1}-\left( \alpha \left( 1-\mu \right) \right)^{n-1}}$. It follows that $d\left( v\right) <v$.
\footnote{To see this precisely, note that the function $\left[ \left( M+v\right)^{\frac{1}{n-1}}\left( \left( n-1\right) M-v\right) \right] -\left[ \left(M-v\right)^{\frac{1}{n-1}}\left( \left( n-1\right) M+v\right) \right] $ is equal to $0$ when $v=0$ and is increasing in $v$ for all $v\geq 0$. Moreover, $\left( M-d\right)^{\frac{1}{n-1}}\left( \left( n-1\right)
M+d\right) $ is decreasing in $d$ for all $d\geq 0$.} Now, $\gamma\left( U\right)$ is convex if and only if $d^{\prime \prime }\left( v\right) \leq 0$. Using (\ref{d def}), $d^{\prime \prime }\left( v\right) \leq 0$ holds if and only if
\[\frac{K^{\prime \prime }\left( \ubar{U}-d\left( v\right)
\right) }{\left( K^{\prime }\left( \ubar{U}-d\left( v\right)
\right) \right)^{2}}\geq \frac{K^{\prime \prime }\left( \ubar{U}+v\right) }{\left( K^{\prime }\left( \ubar{U}+v\right)
\right)^{2}}\text{ .}\]
A direct computation of the derivatives shows that the inequality above holds if and only if
\[\left( \frac{M-d\left( v\right) }{M+v}\right)^{\frac{1}{n-1}}\leq \left( 
\frac{M\left( n-1\right) -d\left( v\right) }{\left( d\left( v\right) \right)
^{2}}\right) /\left( \frac{M\left( n-1\right) +v}{v^{2}}\right) \text{ .}\]
Using (\ref{d def}) again, the inequality above is equivalent to requiring $d\left( v\right) <v$.

It remains to verify that the effective-value density implied by the mixed strategy defined above matches $h_i\left( w\right) $ (given in (\ref{h density})). Observe that the effective values above $\ubar{U}$ are realized as reservation values in the mixed strategy (\ref{mixedstrat42}). For $w\geq \ubar{U}$, the density implied by the mixed strategy is
\[\begin{split}
    &g\left( w\right) \times \frac{\ubar{U}-\gamma\left( w\right) }{w-\gamma\left( w\right) } \\
    &=\frac{\left( 1-\alpha \mu \right)^{n-1}-\left( \alpha \left( 1-\mu\right) \right)^{n-1}}{\left( n-1\right) \hat{U}}\left( \left( \alpha\left(1-\mu \right) \right)^{n-1}+\frac{\left( 1-\alpha \mu \right)^{n-1}-\left( \alpha \left( 1-\mu \right) \right)^{n-1}}{\hat{U}}w\right)^{-\frac{n-2}{n-1}}\frac{w-\gamma\left(w\right) }{\ubar{U}-\gamma\left(w\right) }\times \frac{\ubar{U}-\gamma\left(w\right) }{w-\gamma\left(w\right) } \\
    &= h_{i}\left( w\right)
\end{split} \text{ .}\]
The effective values below $\ubar{U}$ are realized as low posterior realizations in the mixed strategy. Define by $q\colon \left[0,\ubar{U}\right] \rightarrow \left[\ubar{U},\hat{U}\right] $ the inverse of mapping $\gamma$. For $w\leq \ubar{U}$, the density implied by the mixed strategy is
\[\begin{split}
    -q^{\prime }\left( w\right) \times \frac{q\left( w\right) - \ubar{U}}{q\left( w\right) -w}\times g\left( q\left( w\right)\right) &= -\frac{K^{\prime }\left( w\right) }{K^{\prime }\left( q\left( w\right)\right) }\times \frac{q\left( w\right) -\ubar{U}}{q\left(w\right) -w}\times g\left( q\left( w\right) \right) \\
    &=\frac{\left( \left( \alpha \left( 1-\mu \right) \right)^{n-1}+\left(\left( 1-\alpha \mu \right)^{n-1}-\left( \alpha \left( 1-\mu \right)\right)^{n-1}\right) \frac{w}{\hat{U}}\right)^{-\frac{n-2}{n-1}}}{\left(\left(\alpha \left( 1-\mu \right) \right)^{n-1}+\left( \left( 1-\alpha \mu\right)^{n-1}-\left( \alpha \left( 1-\mu \right) \right)^{n-1}\right) \frac{q\left( w\right) }{\hat{U}}\right)^{-\frac{n-2}{n-1}}}\times \frac{\ubar{U}-w}{q\left( w\right) -w}\times g\left( q\left( w\right)\right) \\
    &=h_i\left( w\right)
\end{split} \text{ ,}\]
where the first equality makes use of the definition of the mapping $q$, and
the second and last equalities make use of the definitions of the functions $%
K $ and $g$, respectively.
\end{proof}
\subsection{Proof of Corollary \protect\ref{CS_welfare}}
\begin{proof}
The case of $\mu >\bar{\mu}$ is trivial since the distributions are binary with support $\left\{0,1-c/\mu\right\} $. Consider next the case $\mu <\bar{\mu}$. As shown in Lemma \ref{linear structure no u0}, $\alpha $ is independent of $c$, so it suffices to check the continuous portion of cdf $H$ and $\hat{U}$. The latter is linear and decreasing in $c$, whereas the continuous
portion of $H$ (\ref{linear H no u0}) is obviously increasing in $c$. Finally, according to Corollary 1 of  \cite{choi2}, and appealing to the Law of Iterated Expectations (since we are evaluating welfare from an \textit{ex ante} point of view), the consumer's \textit{ex ante} payoff is given by the expectation of the highest effective value. A worse effective-value distribution in the sense of first-order stochastic dominance thus lowers the consumer's \textit{ex ante} payoff. \end{proof}

\subsection{Proof of Proposition \protect\ref{continuity}}
\begin{proof}
For $c=0$, the effective value is simply the posterior, as any feasible distribution over posteriors has a reservation value equal to one. The symmetric equilibrium distribution $H_{i}^{0}$ of effective values is computed in \cite{Au2} and takes the form
\[H_{i}^{0}\left( w\right) =
\begin{cases}
\left( 1-\alpha_{0}\mu \right) \times \left(\frac{w}{\hat{U}_{0}}\right)^{\frac{1}{n-1}}, \quad &\text{if} \quad w\in \left[ 0,\hat{U}_{0}\right] \\ 
1-\alpha _{0}\mu, \quad &\text{if} \quad w\in \left(\hat{U}_{0},1\right) \\ 
1, \quad &\text{if} \quad w=1
\end{cases}\text{ .}\]
By aligning the slopes of the implied payoff function ($\Pi \left(
1;H_{i}^{0}\right) =\frac{\Pi \left( \hat{U}_{0};H_{i}^{0}\right) }{\hat{U}_{0}}$ if $\alpha _{0}>0$) and the Bayes-plausibility condition ($\int_{0}^{1}wdH_{i}^{0}\left( w\right) =\mu $), the $\alpha _{0}$ and $\hat{U}_0$ can be pinned down as follows. If $\mu \leq \ubar{\mu }=1/n$, $\alpha _{0}=0$ and $\hat{U}_0 = n\mu $. If $\mu >\ubar{\mu }$, then $\hat{U}_{0}=\frac{n\alpha _{0}\mu \left( 1-\alpha _{0}\mu \right)^{n-1}}{1-\left(1-\alpha _{0}\mu \right)^{n}}$, where $\alpha _{0}$ is the unique solution
to $\left( 1-\alpha _{0}\mu \right)^{n}=1-\alpha _{0}$.

For the case $\mu \leq \ubar{\mu }$, the convergence of $H_{i}^{c}$
to $H_{i}^{\ast }$ is immediate by comparing the distribution reported in Lemma \ref{eqdistexplicit} with $H_{i}^{0}$ stated above.

Next, consider the case $\mu \in \left( \ubar{\mu },\bar{\mu}%
\right) $. With $c\rightarrow 0$, $\ubar{U}\rightarrow \mu $;
and hence $H_{i}^{\ast}$ is given by (\ref{linear H no u0}) with $\left(
1-\alpha \mu \right)^{n}-\left( \alpha \left( 1-\mu \right) \right)
^{n}=1-\alpha $, and $\hat{U}=\frac{\left( 1-\alpha \mu \right)
^{n-1}-\left( \alpha \left( 1-\mu \right) \right)^{n-1}}{1/n-\left(
\alpha \left( 1-\mu \right) \right)^{n-1}}\mu $. Recall from the proof of Claim \ref{full info unique} that when $\mu \in \left( \ubar{\mu },\bar{\mu}\right) $, $T\left( \tilde{\alpha}\right) \equiv \left( 1-\tilde{\alpha}\mu \right)^{n}-\left( \tilde{\alpha}\left( 1-\mu \right)\right)^{n}-\left( 1-\tilde{\alpha}\right) $ is negative if and only if $\tilde{\alpha}$ is less than the interior root $\alpha $ of $T$. As $T\left(\alpha _{0}\right) <0$, it follows that $\alpha _{0}<\alpha $. Now, $\Pi\left( w;H_{i}^{0}\right) $ is fully linear over $\left( 0,\hat{U}_{0}\right) $, has a zero vertical intercept, and passes through the points $\left( \mu,1/n\right) $ and $\left( \hat{U}_{0},\left( 1-\alpha _{0}\mu \right)
^{n-1}\right) $; whereas $\Pi \left( w;H_{i}^{\ast}\right) $ is linear over $\left( 0,\hat{U}\right) $, has a positive vertical intercept ($\left( \alpha\left( 1-\mu \right) \right)^{n-1}$), and passes through the points $\left(\mu ,1/n\right) $ and $\left( \hat{U}_{0},\left( 1-\alpha \mu \right)^{n-1}\right) $. Therefore, $\Pi \left( w;H_{i}^{0}\right) $ and $\Pi\left(w;H_{i}^{\ast}\right) $ have a unique intersection at $\mu $ over interior effective values $\left( 0,1\right) $. As $H_{i}^{0}$ and $H_{i}^{\ast}$ have no interior atoms, $\Pi \left( w;H_{i}^{0}\right) =H_{i}^{0}\left( w\right)^{n-1} $and $\Pi \left( w;H_{i}^{\ast}\right) =H_{i}^{\ast}\left( w\right)^{n-1}$ for $w\in \left( 0,1\right) $. As a result, $H_{i}^{0}$ and $H_{i}^{\ast}$ have a unique interior intersection at $\mu $, at which $H_{i}^{0}$ cuts $H_{i}^{\ast}$ from below. As both $H_{i}^{0}$ and $H_{i}^{\ast}$ have the same mean of $\mu $, it follows that $H_{i}^{\ast}$ is a mean-preserving spread of $H_{i}^{0}$.

The case of $\mu \geq \bar{\mu}$ is immediate: $H_{i}^{\ast}$ corresponds to full disclosure, whereas $H_{i}^{0}$ is strictly partial. \end{proof}

\bibliography{sample.bib}

\end{document}